\documentclass[smallextended]{svjour3}       
\smartqed 

\usepackage{graphicx}
\usepackage{subfig}

\usepackage{amsthm}
\usepackage{amsmath}
\allowdisplaybreaks
\usepackage{amsfonts}
\usepackage{amssymb}
\makeatletter
\let\cl@chapter\undefined
\makeatletter
\usepackage{cleveref}
\crefname{equation}{}{}
\usepackage{enumitem}
\usepackage{longtable}
\usepackage{booktabs}
\usepackage{nomencl}
\usepackage{etoolbox}
\patchcmd{\thenomenclature}{\section*}{\section}{}{}
\makenomenclature

\usepackage{xstring}
\usepackage{xpatch}

\usepackage{marginnote}
\usepackage{lineno}
\newcommand{\cR}{\mathbb{R}}
\newtheorem{assum}[theorem]{Assumption}
\newcommand{\eps}{\epsilon}
\newcommand{\be}{\begin{equation}}
\newcommand{\ee}{\end{equation}}

\newcommand{\tpk}{{p_k}}

\newcommand{\grad}{\nabla}

\newcommand{\ve}{{\gamma}}

\usepackage[algo2e,ruled,linesnumbered,vlined]{algorithm2e}

\newcommand{\lh}{{\hat l}}

\begin{document}

\title{A Trust Region Method for the Optimization of Noisy Functions\thanks{Sun was supported by NSF grant DMS-1620022. Nocedal was supported by AFOSR grant FA95502110084 and ONR grant N00014-21-1-2675.}
}
\author{Shigeng Sun     \and
    Jorge Nocedal 
}
\institute{Shigeng Sun \at
       Department of Engineering Sciences and Applied Mathematics, Northwestern University, Evanston, IL, USA \\
       \email{shigengsun2024@u.northwestern.edu}      
      \and
      Corresponding author: Jorge Nocedal \at
       Department of Industrial Engineering and Management Sciences, Northwestern University, Evanston, IL, USA\\
       \email{j-nocedal@northwestern.edu}
}
\date{Received: 02 January 2022 / Accepted: }
\maketitle
\begin{abstract}
Classical trust region methods were designed to solve problems in which function and gradient information are exact. This paper considers the case when there are bounded errors (or noise) in the above computations and proposes a simple modification of the trust region method to cope with these errors. The new algorithm only requires information about the size of the errors in the function evaluations and incurs no additional computational expense. It is shown that, when applied to a smooth (but not necessarily convex) objective function, the iterates of the algorithm visit a neighborhood of stationarity infinitely often, and that the rest of the sequence cannot stray too far away, as measured by function values. Numerical results illustrate how the classical trust region algorithm may fail in the presence of noise, and how the proposed algorithm ensures steady progress towards stationarity in these cases.
\end{abstract}

\section{Introduction}
Trust region methods are powerful techniques for nonlinear optimization that have the ability to incorporate second-order information, without requiring it to be positive definite. They are endowed with strong global convergence properties and have proven to be effective in practice. Although the design and analysis of trust region methods are well established in the absence of noise (or errors), this is not the case when noise is present.

In this paper, we show how to redesign the classical trust region method for unconstrained optimization to handle problems where the objective function, gradient, and (possibly) Hessian, are subject to bounded, non-diminishing noise. This involves only one modification in the algorithm: the ratio of actual/predicted reduction used for step acceptance is now relaxed by a term proportional to the noise level. All other aspects of the classical trust region method remain unchanged. We show that, under mild conditions, the proposed algorithm converges to a neighborhood of stationary points, where the size of the neighborhood is determined by the level of noise. This analysis is more complex than for line search methods due to the effects of memory encapsulated in the trust region update. Our convergence results do not assume convexity of the objective function but only that it is sufficiently smooth.

Examples of practical optimization applications with bounded noise include those that employ mixed-precision arithmetic; problems where derivatives are approximated by finite differences; and problems in which the evaluation of the objective function (and gradient) contain computational noise. 

This investigation was motivated by numerical experiments performed by the authors that indicated that, although the classical trust region approach often tolerates significant levels of noise, it can fail in certain situations. This raises the question of how to best modify the method to avoid failures. The algorithm proposed here is inspired by work on line search methods for unconstrained optimization \cite{berahas2019derivative,xie2020analysis} and equality constrained optimization \cite{figen2021ConDFO}. In those papers, convergence-to-neighborhood results were derived but the analysis presented here follows different lines, as trust region methods require different proof techniques.

The paper is organized into 5 sections. In the rest of this section, we provide a review of the relevant literature. In section 2, we describe the problem setting and the proposed trust region algorithm. The main convergence results are presented in section~3. Numerical experiments, summarized in section~4, indicate that the new algorithm {is more robust than} the classical method. Section~5 presents the final remarks on the contributions of this work.

\subsection{Literature Review}
The study of nonlinear optimization problems with errors or noise in the function and gradient has attracted attention in recent years, motivated by 
 the use of finite difference approximations to derivatives \cite{nesterov2017random,shi2021numerical,shi2021automatic} and by applications in machine learning; see \cite{curtis2020adastoc} for a review of some recent work. 

One of the earliest investigations of trust region methods with errors is \cite{carter}, which proved global convergence assuming that the errors in the gradient diminish at a rate that is proportional to the norm of the true gradient; this condition is referred to as the \emph{norm test} in \cite{bollapragada2018adaptive,exact2018ima}. The importance of the norm test was promoted in \cite{byrd2012sample}, which established linear convergence and complexity bounds for an adaptive sampling method for empirical risk minimization, as well as in \cite{cartis2018global,paquette}, which establishes convergence in probability for a stochastic line search method.  

Prior studies of optimization methods for minimization of functions with non-diminishing, bounded errors include \cite{berahas2019derivative}, which employed a relaxed Armijo back-tracking line search and established linear convergence to a neighborhood of the solution for strongly convex functions. Stopping time guarantees for the same relaxed line search is proven in \cite{berahas2019global}.
A similar relaxed Armijo back-tracking line search technique is used in \cite{katya2021stolinesearch}, which considered different oracles from \cite{paquette} to allow biased estimates, and provided complexity bounds for different noise structures under probabilistic frameworks. Quasi-Newton methods were analyzed in \cite{xie2020analysis}, which described a noise tolerant modification of the BFGS method; \cite{shi2020noise} showed ways to make this method robust and efficient in practice. 

For constrained optimization, \cite{berahas2021stochastic,berahas2021sequential,curtis2021inexact} studied a sequential quadratic programming (SQP) method for equality constrained optimization in the case when the objective function is stochastic and the constraints are deterministic. Those three papers give conditions under which convergence can be expected, giving careful attention to the behavior of the penalty parameter. 
Using a relaxed Armijo line search procedure, \cite{figen2021ConDFO} shows global convergence to a neighborhood of the solution for an SQP method for equality constrained problems.

Analysis for trust region methods with more general (unbounded) noise is presented in \cite{katya2018storm}, which establishes almost sure global convergence under the assumption that function and gradient information is sufficiently accurate with high enough probability. \cite{cartis2019TRsupermartingale} views the optimization as a generic stochastic process, and improves upon the results of \cite{katya2018storm}. The analysis presented in \cite{cartis2019TRsupermartingale}
establishes convergence results for a trust region method and, under the assumption of sufficiently accurate stochastic gradient information, derives a stopping time result and a second order global complexity bound. A method inspired by trust region techniques is \cite{curtis2019stochTR}, which uses step normalization techniques in the stochastic optimization setting, and establishes conditions for linear and sublinear convergence.
A series of papers, including \cite{cartis2021TRqDA,toint2021TRqEDAN,toint2021ARqpEDA2}, analyze regularization and {trust region methods} with adaptive accuracy in the function and gradient evaluations, and establish worst case complexity bounds. 

The style of analysis presented in \cite{katya2018storm,cartis2019TRsupermartingale,curtis2019stochTR}, 
which is used to prove convergence in probability,  
stands in contrast with the deterministic technique employed in this paper, which assumes bounded errors. It remains to be seen which approach is more useful for the design of noise tolerant optimization methods---or whether the two approaches complement each other.

\section{ Problem Statement and Algorithm}
\label{sec:Algorithm}
Our goal is to design a trust region method to solve the unconstrained minimization problem
\be
   \min_{x\in\cR^n} f(x), \label{eq:tr-prob}
\ee
in the case when the function $f(x)$ and gradient $g(x) = \grad f(x)$
 cannot be evaluated exactly. Instead, we have access to noisy observations of the above quantities, which we denote as $\tilde f(x)$, and $\tilde g(x)$. 
 We write
\be
\tilde f(x)  = f(x) + \delta_f(x),\quad\text{and}\quad
\tilde g(x) = g(x) + \delta_g(x), \label{eq:tr-noise}
\ee
where the error functions (or noise) $ \delta_f(x)$, $\delta_g(x)$
 are assumed to be bounded, i.e., 
 \be\label{eq:tr-asp1}
|\delta_f(x)|\leq \eps_f, \qquad
\|\delta_g(x)\|\leq \eps_g, \qquad \ \forall x\in\cR^n.
\ee
Throughout the paper $\| \cdot \|$ stands for the Euclidean norm.

Let us apply a classical trust region method to problem \cref{eq:tr-prob}. At each iterate, the method constructs a quadratic model
\be \label{eq:tr-nmodel} 
   m_{k}(p)=\tilde f(x_k)+\tilde g(x_k)^{T} p+\frac{1}{2} p^{T} \tilde B_{k} p ,  
\ee
and solves the following trust region subproblem for the step $p_k$:
\be \min _{p \in \mathbb{R}^{n}} m_{k}(p)\quad \text { s.t. }\|p\| \leq \Delta_{k}\label{eq:tr-subprob}. \ee
In \Cref{eq:tr-nmodel}, $\tilde B_k$ could be defined as a noisy evaluation of the Hessian, a quasi-Newton matrix, or some other approximation. 
To decide if the step $p_k$ should be accepted---and if the trust region radius $\Delta_k$ should be modified--- classical trust region methods employ the ratio of actual to predicted reduction in the objective function, defined as
\be \label{eq:old}
 \frac{\tilde f\left(x_{k}\right)-\tilde f\left(x_{k}+p_{k}\right)}{m_{k}(0)-m_{k}\left(p_{k}\right)}. 
\ee
This ratio is, however, not adequate in the presence of noise because if $\Delta_k$ becomes very small, the numerator can be of order $\eps_f$, while the denominator will be proportional to $\Delta_k$. Thus, if $\Delta_k << \epsilon_f$, the ratio \cref{eq:old} may exhibit wild oscillations that can cause the algorithm to perform erratically; see the examples in Section~\ref{sec:numerical}.

To address this issue, we propose the following noise tolerant variant of \cref{eq:old}:
\be
  \rho_{k}=\frac{\tilde f\left(x_{k}\right)-\tilde f\left(x_{k}+p_{k}\right)+r\eps_f}   
  {m_{k}(0)-m_{k}\left(p_{k}\right)+r\eps_f} \label{eq:tr-ratio},
\ee
where $r > 2 $ is a constant specified below.
The reason for relaxing both the numerator and denominator in \cref{eq:tr-ratio} is to be consistent with the classical narrative of trust region methods where a ratio close to 1 is an indication that the model is adequate. An alternative approach would be to relax only the numerator and interpret the condition $\rho_k > c$ (where $c>0$ is a constant) as a relaxed Armijo condition of the type studied in \cite{berahas2019derivative,figen2021ConDFO}. We find the first interpretation to be easier to motivate and to yield tighter bounds in the convergence analysis. We state the algorithm as follows.

\bigskip

\begin{algorithm2e}[H] \label{algorithm}
\SetAlgoLined
 Initialize $\Delta_{0}$, and chose constants $ 0<c_0\leq c_1<c_2<1$ and $\nu > 1$
 
 \While{a termination condition is not met}{
 Compute $p_{k}$ by solving \cref{eq:tr-subprob} (exactly or approximately);\\
 Evaluate $\rho_{k}$ as in \cref{eq:tr-ratio}; \\
  
 \uIf{$\rho_{k}<c_1$}{
  $\Delta_{k+1} = \frac1\nu\Delta_{k}$;
 }
  \uElseIf{$\rho_{k}>c_2$ }{
  $\Delta_{k+1}= \nu \Delta_{k}$;
 }
 \Else{
  $\Delta_{k+1}=\Delta_{k}$;
 }  
 \uIf{$\rho_{k}>c_0$}{
  $x_{k+1} = x_k+p_k $;
 }\Else{
  $x_{k+1} = x_k$; 
 }
 Set $k \leftarrow k+1$;
 }
\caption{Noisy Trust-Region Algorithm}
\end{algorithm2e}

\bigskip

Typical values of the parameters are $c_0= 0.1$, $c_1=\frac14$, $c_2=  \frac12$, $\nu = 2$, but other values can be used in practice. The global convergence result presented in the next section holds if the constant $r$ in \cref{eq:tr-ratio} is chosen as
\be r = 2/{(1-c_2)}\label{eq:rdef}.\ee
We assume that the step $p_k$ computed in step 3 yields a decrease in the model $m_k$ that is at least as large as that given by the Cauchy step (defined below). This provides much freedom in the design of the algorithm, and includes the dogleg and Newton-CG methods, as well as the exact solution of the trust region problem; see, e.g., \cite{mybook}. 

In practice it can be useful to increase the trust region radius in Step~7 only if $\rho_{k}>c_2$ and $\|p_{k}\| = \Delta_{k}$, as this can prevent unnecessary oscillations in the trust region radius. 
The convergence result presented in the next section can {easily} be extended to that case, assuming certain technical conditions on the step computation---which are satisfied by the dogleg and Newton-CG methods. 

\section{Global Convergence Analysis} 
\label{sec:convergence}
In this section, we establish a global convergence result for Algorithm~\ref{algorithm} that applies to general objective functions. The proof is based on the observation that, when the gradient is large enough, the trust region radius will eventually become large too, ensuring sufficient descent in the objective function despite the presence of noise. This drives the iteration toward regions where the stationarity measure is small (i.e., comparable to the noise level). 

We begin by establishing a standard requirement on the step computation based on the \emph{Cauchy step} $p_{k}^{c}$ for problem \cref{eq:tr-prob}, which is defined as
\be 
   p_{k}^{c}=- \tau_{k} \frac{\Delta_{k}}{\left\|\tilde g_{k}\right\|} \tilde g_{k},\label{eq:cauchy}
\ee
where
\be \label{eq:cauchy1}
\tau_{k}=\left\{\begin{array}{ll}
1 & \text { if } \tilde g_{k}^{T} \tilde B_{k} \tilde g_{k} \leq 0 \\
\min \left(\left\|\tilde g_{k}\right\|^{3} \big/\left(\Delta_{k} \tilde g_{k}^{T} \tilde B_{k} \tilde g_{k}\right), 1\right) & \text { otherwise.}
\end{array}\right.
\ee
As is well known (see e.g. \cite[Lemma 4.3]{mybook}), the reduction in the model provided by the Cauchy step satisfies 
\be m_k(0) - m_k(p_k^c)\geq \frac{1}{2}\left\|\tilde g_{k}\right\| \min \left(\Delta_{k}, \frac{\left\|\tilde g_{k}\right\|}{\left\|\tilde B_{k}\right\|}\right).\label{cdecrease}\ee
We assume that the step $p_k$ computed by Algorithm~\ref{algorithm} yields a reduction in the model that is not less than that produced by the Cauchy step, i.e.,
\be  \label{eq:mred}
 m_k(0) - m_k(p_k)\geq m_k(0) - m_k(p_k^c)\geq \frac{1}{2}\left\|\tilde g_{k}
 \right\| \min \left(\Delta_{k}, \frac{\left\|\tilde g_{k}\right\|}{\left\|\tilde B_{k}
 \right\|}\right).
\ee

We can now state the assumptions on the problem and the algorithm under which the global convergence results are established.
\begin{assum}\label{assum1}
The objective function $f$ is Lipschitz continuously differentiable with constant $L$, i.e.,
\be\label{eq:lip}
  \| g(x) - g(y)\| < L \|x-y\|. 
\ee 
\end{assum}
\begin{assum}\label{assum2}
\noindent The error in the function and gradient evaluations is bounded, i.e., \cref{eq:tr-asp1} holds for some constants $\eps_f, \eps_g.$
\end{assum}
\noindent We impose no other conditions on the errors, other than boundedness. Next, we impose a minimal requirement on the Hessian approximations.
\begin{assum}\label{assum3}
There is a constant $L_B >0$ such that the matrices $\tilde B_k$ satisfy
\be \|\tilde B_k\| < L_B , \ \forall k . \label{eq:lb}\ee
\end{assum}
\noindent There is freedom in the computation of the step $p_k$, but it must yield Cauchy decrease. 

\begin{assum}\label{assum4}
The step $p_k$ computed by Algorithm~\ref{algorithm} satisfies (\ref{eq:mred}).
\end{assum}
\noindent {This assumption can be relaxed so as to require only a fraction of Cauchy decrease, but we do not do so here to avoid the introduction of more constants.}
The final requirement is standard.
\begin{assum}\label{assum5}
The sequence $\{ \tilde f_k\}$ generated by Algorithm~\ref{algorithm} is bounded below.
\end{assum}

 \noindent We now proceed with the analysis.

\subsection{Properties of the ratio $\rho_k$}
We begin by establishing a bound between $\rho_k$ and 1. From \cref{eq:tr-ratio}, we have
\be 
\begin{aligned}
\label{eq:omegaabs}
\left|\rho_k-1\right|&=\left|\frac{m_{k}\left(p_{k}\right)-\tilde f\left(x_{k}+p_{k}\right)}{m_{k}(0)-m_{k}\left(p_{k}\right)+r\eps_f}\right| .
\end{aligned}
\ee
From Taylor's Theorem we have
\begin{align*}
\tilde f(x_k+p_k) & =f(x_k+p_k) + \delta_f(x_k +p_k) \\
  & = f\left(x_{k}\right)+g_k^{T} p_{k}+\int_{0}^{1}\left[g\left(x_{k}+t p_{k}\right)-g_k\right]^{T} p_{k} d t +\delta_f(x_k+p_k).
 \end{align*}
With this, by \cref{eq:lip}, \cref{eq:lb}, and \cref{eq:tr-asp1}, we obtain
\begin{align} \label{eq:tr-top}
\left|m_{k}(p_k)-\tilde f(x_k+p_k)\right| 
\leq & \ \tfrac{1}{2} (L_B+L)\|\tpk\|^2 + \eps_g\|\tpk\| + 2\eps_f\\
\equiv & \ M\|\tpk\|^2 + \eps_g\|\tpk\| +2\eps_f , \nonumber
\end{align}
where 
\be M = \tfrac{1}{2} (L_B+L). \label{eq:tr-nu}\ee
By substituting \cref{eq:tr-top} and \cref{eq:mred} into \cref{eq:omegaabs}, we establish the following result.

\begin{lemma}\label{tr-modelred6}
If $\rho_k$ is defined by \cref{eq:tr-ratio}, then for all $k$,
\be
\left|\rho_k-1\right| \leq \frac{ M \Delta_k^2 + \eps_g \Delta_k+ 2\eps_f}{\frac12\|\tilde g_k\|\min(\Delta_k,\|\tilde g_k\|/\|\tilde B_k\|)+r\eps_f}.\label{eq:tr-rhodist}
\ee
\end{lemma}

This lemma suggests that $\rho_k$ can be made close to 1 by decreasing $\Delta_k$, up until the noise term $\eps_f$ dominates. This assertion will be made more precise below.

\subsection{Lower Bound on Trust Region Radius}
We now show that if $\Delta_k$ is very small and the gradient is large compared to the noise $\eps_g$, Algorithm~\ref{algorithm} will increase the trust region radius. We recall that $r$ is defined in \cref{eq:rdef} and that $\nu > 1$.

\begin{lemma}[Increase of Trust Region Radius]\label{tr-lb}
Suppose that, at iteration $k$,
\be \| \tilde g_k \|> r\eps_g+\ve , \label{eq:gbig}\ee
 for some constant $\ve>0$. Then, if 
\be \label{eq:dhat} 
 \Delta_k \leq \bar\Delta =: \frac{\ve}{r M },
 \ee
 we have that
\be \Delta_{k+1} = \nu\Delta_k . \ee
\end{lemma}

\begin{proof} 
Since $r > 2$, we have from \cref{eq:lb}, \cref{eq:tr-nu} and \cref{eq:gbig} that
\be r M > 2 M >\|\tilde B_k\| \quad\text{and }\quad \ve < \|\tilde g_k\|, \ee and thus
\be \label{eq:can} 
  \bar\Delta< \|\tilde g_k\|/\|\tilde B_k\|.
\ee

Thus, if $\Delta_k\leq\bar\Delta$, we have
\be \min(\Delta_k ,\|\tilde g_k\|/\|\tilde B_k\|) = \Delta_k\label{eq:tr-compare}.\ee
In addition, if $\Delta_k\leq\bar\Delta$, we also have
\be
 M \Delta_k+ \eps_g\leq M \bar\Delta+ \eps_g = \frac{\ve}{r}+ \eps_g =\frac1r(r\eps_g+\ve). \label{eq:tr-quant1}
\ee
Substituting \cref{eq:tr-compare}, \cref{eq:gbig}, \cref{eq:tr-quant1} and \cref{eq:rdef} into \cref{eq:tr-rhodist}, we have that for all $\Delta_k\leq \bar\Delta$
\begin{align}
\left|\rho_{k}-1\right| &\leq \frac{ M \Delta_k^2 + \eps_g \Delta_k+ 2\eps_f}{\frac12\|\tilde g_k \|\Delta_k+r\eps_f}\nonumber\\
& < \frac{ M \Delta^2_k + \eps_g \Delta_k+ 2\eps_f}{\frac12(r\eps_g+\ve)\Delta_k+r\eps_f}\nonumber\\
&< \frac{\frac1r(r\eps_g+\ve) \Delta_k+ 2\eps_f}{\frac12(r\eps_g+\ve)\Delta_k+r\eps_f}\nonumber\\
& = \frac2r\nonumber\\
& = 1-c_2.\label{eq:tr-quant2}
\end{align}
This implies that $\rho_k > c_2$, and by step~8 of Algorithm~\ref{algorithm} we have that $\Delta_{k+1}= \nu \Delta_k.~~~$
\end{proof}

A consequence of this lemma is that there is a lower bound for the trust region radius if the norm of the noisy gradient remains greater than $r\eps_g$. 

\begin{corollary}[Lower Bound on Trust Region Radius]\label{tr-lbprop}
Given $\ve>0$, if there exist $K>0$ such that for all $k\geq K$
\be \| \tilde g_k \|> r\eps_g+\ve,\label{eq:tr-ctasp}\ee
then there exist $K_0 \geq K $ such that for all $k\geq K_0$,
\be \Delta_k > \tfrac{1}{\nu}\bar\Delta= \frac{\ve}{\nu r M }.\label{eq:tr-bardel} \ee
\end{corollary}
\begin{proof}
We apply \Cref{tr-lb} for each iterate after $K$ to deduce that, whenever $\Delta_k \leq \, \bar \Delta$, the trust region radius will be increased. Thus, there is an index $K_0$ for which $\Delta_k$ becomes greater than $\bar\Delta$. On subsequent iterates, the trust region radius can never be reduced below $\bar\Delta/\nu$ (by Step~6 of Algorithm~\ref{algorithm}) establishing the bound \cref{eq:tr-bardel}. 
\end{proof}

{\em Remark.} In traditional trust region analysis for deterministic (noiseless) optimization, one shows that the trust region radius will not shrink below a certain value that depends on the Lipschitz constant and the norm of the current gradient. However, that analysis does not imply that the trust region will increase beyond a certain threshold, which is required in the presence of noise. We need to show that the trust region eventually becomes large enough with respect to the noise level so that progress can be made. This differentiates our analysis from classical trust region convergence theory.

\subsection{Reduction of Noisy Function}
The classical trust region algorithm is monotonic, as it requires a reduction in the objective function when accepting a step. Due to the relaxation in \cref{eq:tr-ratio}, Algorithm~\ref{algorithm} can accept steps that increase the noisy function. However, when the iterates are far from the solution, this is not the case. We now show that when the noisy gradient and trust region radius are both large enough, the reduction in the objective is large enough to overcome any increase allowed by \cref{eq:tr-ratio}.

\begin{lemma}[Noisy Function Reduction]\label{noisyfuncred}
Suppose that for some $k>0$ 
\be \label{eq:supp}
   \| \tilde g_k \|> r\eps_g+ \gamma \quad\mbox{and}\quad
   \Delta_k \, {\geq} \, \frac{\bar\Delta}{\nu} = \frac{\gamma}{\nu rM} ,
 \ee 
 where 
 \begin{equation} \label{eq:gammap}
 \gamma= \eta+ \mu ,
 \end{equation} 
  with $\mu >0$ an arbitrarily small constant, and 
 \be \label{eq:deldef}
 \eta = \frac12\left(- r\eps_g+\beta \right),\quad \beta = \sqrt{(r\eps_g)^2 + 8\nu r^2\left(\frac1{c_0}-1\right) M \eps_f }.
 \ee
Then, if the step is accepted at iteration $k$ {by Algorithm~\ref{algorithm}}, we have 
\be  
\tilde f\left(x_{k}\right)-\tilde f\left(x_{k}+p_{k}\right) > \frac{c_0}{2\nu r M }\left(\mu\beta+ \mu^2\right).
\ee
\end{lemma}

\begin{proof} 
As argued in \cref{eq:can}, $\bar \Delta= \frac{\ve}{rM}< \frac{\left\|\tilde g_{k}\right\|}{\left\|\tilde B_{k}\right\|} $, and therefore
\be\label{eq:tr-lem36_1}\min \left(\Delta_{k}, \frac{\left\|\tilde g_{k}\right\|}{\left\|\tilde B_{k}\right\|}\right) \geq \frac{\ve}{\nu rM}.\ee
If the step $p_k$ is accepted, we have from Step~12 of Algorithm~\ref{algorithm} that $\rho_k>c_0$, which by \cref{eq:tr-ratio} is equivalent to
\be \frac{\tilde f\left(x_{k}\right)-\tilde f\left(x_{k}+p_{k}\right)+r\eps_f}{m_{k}(0)-m_{k}\left(p_{k}\right)+r\eps_f}>c_0 .\ee
Thus by \cref{eq:mred}, \cref{eq:supp}, \cref{eq:tr-lem36_1} and \cref{eq:gammap}
\begin{align}\tilde f\left(x_{k}\right)-\tilde f\left(x_{k}+p_{k}\right)
 > &c_0\left[m_k(0) -m_k(p_k)\right] + r(c_0-1)\eps_f \nonumber\\
\geq&\frac{{c_0}}{2}\left\|\tilde g_{k}\right\| \min \left(\Delta_{k}, \frac{\left\|\tilde g_{k}\right\|}{\left\|\tilde B_{k}\right\|}\right)+ r(c_0-1)\eps_f\nonumber\\
> &\frac{{c_0}}{2\nu r M }\left(r\eps_g + \ve\right)\ve+ r(c_0-1)\eps_f\nonumber \\
> &\frac{{c_0}}{2\nu r M }\left(r\eps_g + \eta\right)\eta+ r(c_0-1)\eps_f .
\label{eq:tr-nred1}
\end{align}
{We now chose $\eta$ so that the right hand side is positive. We obtain}
%
\[ 
 \eta \geq \frac12\left(- r\eps_g+\beta \right) \quad
 \text{or} \quad
 \eta \leq \frac12\left(- r\eps_g-\beta \right) 
 \]
 {We wish for $\eta$ to be the smallest positive value satisfying these inequalities, yielding}
\be \eta = \frac12\left(- r\eps_g+\beta \right). \ee
{Substituting this quantity in} \cref{eq:tr-nred1}, we have
\begin{align*}\tilde f\left(x_{k}\right)-\tilde f\left(x_{k}+p_{k}\right)
> &\frac{{c_0}}{2\nu r M }\left(r\eps_g + \ve\right)\ve+ r(c_0-1)\eps_f\\
= &\frac{{c_0}}{2\nu r M }\left(r\eps_g + \eta + \mu \right)(\eta + \mu)+ r(c_0-1)\eps_f\\
=&\frac{{c_0}}{2\nu r M }\left(r\eps_g + \frac12\left(- r\eps_g+\beta \right)+ \mu \right)\left(\frac12\left(- r\eps_g+\beta \right) + \mu\right)+ r(c_0-1)\eps_f\\
=&\frac{{c_0}}{2\nu r M }\left(r\eps_g/2 + \beta/2+ \mu \right)\left(- r\eps_g/2+\beta/2 + \mu\right)+ r(c_0-1)\eps_f\\
=&\frac{{c_0}}{2\nu r M }\left[\left(\beta/2+ \mu \right)^2 - \left( r\eps_g/2\right)^2\right]+ r(c_0-1)\eps_f\\
=&\frac{{c_0}}{2\nu r M }\left[(\beta/2)^2+\mu\beta +\mu ^2 - \left( r\eps_g/2\right)^2\right]+ r(c_0-1)\eps_f\\
=&\frac{{c_0}}{2\nu r M }\left[\frac{\beta^2- (r\eps_g)^2}{4}+\mu\beta +\mu ^2 \right]+ r(c_0-1)\eps_f\\
=&\frac{{c_0}}{2\nu r M }\left[\frac{(r\eps_g)^2 + 8\nu r^2\left(\frac1{c_0}-1\right) M \eps_f - (r\eps_g)^2}{4}+\mu\beta +\mu ^2 \right]+ r(c_0-1)\eps_f\\
=&\frac{{c_0}}{2\nu r M }\left[{ 2\nu r^2\left(\frac1{c_0}-1\right) M \eps_f}+\mu\beta +\mu ^2 \right]+ r(c_0-1)\eps_f\\
=& r(1-c_0)\eps_f + \frac{{c_0}}{2\nu r M }\left(\mu\beta +\mu ^2 \right)+ r(c_0-1)\eps_f\\
=& \frac{c_0}{2\nu r M }\left(\mu\beta+ \mu^2\right).
\end{align*}

\end{proof}
{The first inequality \cref{eq:supp}, together with \cref{eq:gammap}, \cref{eq:deldef}, identify the region where noise does not dominate and progress in the objective function can be guaranteed.
The constant $\mu$ was introduced to ensure that our analysis is meaningful in the case when noise is not present ($\epsilon_f=\epsilon_g=0$), as it shows that a decrease in the objective is achieved. Nonetheless, the global convergence results presented below are of interest only when noise is present, so there we essentially absorb $\mu$ into $\eta$ by setting $\mu = \epsilon_g/2$.}

To summarize the results obtained so far, \Cref{tr-lb} states that when $\|\tilde g_k\|$ is large enough, the trust region is either large enough or will eventually be increased to be so.  \Cref{noisyfuncred} states that when the gradient and trust region are both large enough, every accepted iterate reduces the noisy objective function by a non-vanishing amount. We show that this drives iterations towards stationary points of the problem.

\subsection{{Global Convergence Theorems}}
Our global convergence results are presented in two parts. The first result states that the iterates visit, infinitely often, a critical region characterized by a small gradient norm. The second result states that after visiting the above critical region for the first time, the iterates cannot stray too far from it, as measured by the objective value.

\begin{theorem}[Global Convergence to Critical Region]\label{thm:GCT1}
Suppose that \Cref{assum1} through \Cref{assum5} are satisfied.
Then, the sequence of iterates $\{x_k\}$ generated by Algorithm~\ref{algorithm} visits infinitely often the critical region ${C}_1$ defined as
\be {C}_1 =\left\{x: \|g(x)\| \leq \left(r+1\right) \eps_g + \frac{ \beta}{2} \right\} ,
\label{eq:C1}
\ee
{where $r$ and $\beta$ are defined in \cref{eq:rdef}, \cref{eq:gammap}, \cref{eq:deldef}, with $\mu= \epsilon_g/2$, $\nu > 1$ and $M$ given by \cref{eq:tr-nu}.}
\end{theorem}

\begin{proof}
Assume by way of contradiction that there exist $K'$ such that for all $k>K'$
\be \|g(x_k)\| > \left(r+1\right) \eps_g + \frac\beta2\label{eq:palm}.\ee
Thus, {by \cref{eq:tr-asp1}}, definition \cref{eq:deldef} of $\eta$, and setting $\mu = \eps_g/2$, we have that for all $k > K'$
\begin{align}
 \|\tilde g(x_k)\| & > r \eps_g + \tfrac{1}{2}\beta\nonumber \\
 	& = -\tfrac{1}{2}r\eps_g+ \tfrac{1}{2}\beta + \tfrac{3}{2}{r\eps_g}\nonumber\\
	& = \eta + r\eps_g + \tfrac{1}{2}{r\eps_g}\nonumber \\
	& > \, r\eps_g + \eta + \mu \qquad\qquad\mbox{(since $r>1$)}\nonumber \\
 	& = r\eps_g+\ve. \qquad\quad\qquad\mbox{~~(by \cref{eq:gammap})} \label{eq:tr-g1}
   \end{align}
 
We now apply \Cref{tr-lbprop} and deduce that there exist $K_0\geq K'$, such that for all $k \geq K_0$,
\be \Delta_k > \frac{\ve}{\nu r M }. \label{eq:furn} \ee   
    
When a step is not accepted, $\rho_k < c_0 < c_1$, and Algorithm~\ref{algorithm} will reduce the trust region radius. If no step is accepted for all $k > K_0$, the trust region radius would shrink to zero, contradicting \cref{eq:furn}. Therefore, there must exist infinitely many accepted steps. Now, by \cref{eq:tr-g1}, \cref{eq:furn} the conditions of \Cref{noisyfuncred} hold, and we deduce that each accepted step  $ k'>K_0$ achieves the reduction
\be  
\tilde f\left(x_{k'}\right)-\tilde f\left(x_{k'}+p_{k'}\right)>
 \frac{c_0}{2\nu r M }\left(\mu\beta+ \mu^2\right) = \frac{c_0}{2\nu r M }\left(\frac{\eps_g}{2}\beta+\frac{\eps_g^2}{4}\right).
\ee
Since, as mentioned above, there is an infinite number of accepted steps, we deduce that $ \{\tilde f(x_k)\} \rightarrow - \infty $, contradicting Assumption~\ref{assum5}. Therefore, the index $K'$ defined above cannot exist and we have that \cref{eq:palm} is violated an infinite number of times. 
\end{proof}

The achievable accuracy in the gradient guaranteed in \cref{eq:C1} depends on $\eps_g$ and $\sqrt{\eps_f}$, by the definition of $\beta$. The dependence on $\eps_g$ is evident, while the dependence on $\sqrt{\eps_f}$ is due to the combined (multiplicative) effect of the gradient and the trust region radius bound.

{Before stating our next theorem, we prove two simple technical results.}

\begin{proposition}\label{acceptedstep}
If Algorithm~\ref{algorithm} takes a (nonzero) step at iteration $k$, then
\be \tilde f_{k+1} - \tilde f_k < r(1-c_0)\eps_f .\ee
\end{proposition}
\begin{proof}

If the step is taken, we have from Step~12 of Algorithm~\ref{algorithm} that $\rho_k>c_0$, which by \cref{eq:tr-ratio} is equivalent to
\be \frac{\tilde f\left(x_{k}\right)-\tilde f\left(x_{k}+p_{k}\right)+r\eps_f}{m_{k}(0)-m_{k}\left(p_{k}\right)+r\eps_f}>c_0 ,\ee
and since $p_k$ cannot increase the model $m_k$, we have
\be \tilde f\left(x_{k}\right)-\tilde f\left(x_{k}+p_{k}\right) > c_0\left[m_k(0) -m_k(p_k)\right] + r(c_0-1)\eps_f > r(c_0-1)\eps_f.\ee

\end{proof}

{Next, we employ \Cref{tr-lb} and obtain the following result.}
\begin{corollary}[Maintaining Lower Bound on Trust Region Radius]\label{tr-lbcor2}
Let $\ve>0$ be defined by \cref{eq:gammap}--\cref{eq:deldef}, and suppose there exist $K>0$ and $\hat K > K$ such that for $k = K+1, ..., \hat K-1$
\be \| \tilde g_k \|> r\eps_g+\ve , \ee
and that 
\be \Delta_{K+1} \geq \frac{\gamma}{\nu rM} = \frac{\bar \Delta}{\nu}.\ee
Then for $k = K+1, ..., \hat K-1$
\be \label{eq:maint} 
\Delta_k \geq \frac{\gamma}{\nu rM}=\frac{\bar \Delta}{\nu}. \ee
\end{corollary}

\begin{proof}
The proof is by induction. Condition \cref{eq:maint} holds for $k= K+1$. We show that if \cref{eq:maint} it holds for some $k \in \{K+1,\ldots, \hat K-2\}$, then it holds for $k+1$.

Specifically, suppose that for such $k$ we have that
\be \Delta_k \geq \frac{\gamma}{\nu rM}\label{eq:basestep}. \ee
By \Cref{tr-lb}, if $\Delta_k \leq \frac{\gamma}{rM}$, the trust region radius is increased, i.e.,  
\be \Delta_{k+1} = \nu\Delta_k \geq \frac{\gamma}{ rM}>\frac{\gamma}{\nu rM}\label{eq:induction1} .\ee 
 If on the other hand $\Delta_k > \frac{\gamma}{rM}$, the trust region radius could be decreased, but in that case 
\be \Delta_{k+1} \geq \frac{\Delta_k}{\nu}>\frac{\gamma}{\nu rM}.\label{eq:induction2}\ee 
\end{proof}

The next theorem shows that after an iterate has entered the neighborhood $C_1$ defined in \Cref{thm:GCT1}, all subsequent iterates cannot stray too far away in the sense that their function values remain within a band of the largest function value in $C_1$.

\begin{theorem}[Iterates Remain in the Level Set $C_2$]\label{thm:GCT2i}
Suppose that \Cref{assum1} through \Cref{assum5} are satisfied. Then, after the iterates $x_k$ generated by Algorithm ~\ref{algorithm} visit $C_1$ for the first time, they never leave the set $C_2$ defined as
\be \label{eq:c2def} {C}_2 =\left\{x: f(x) \leq\sup_{y\in C_1}  f(y) + 2\eps_f+ \max[G, r(1-c_0)\eps_f ]\right\},\ee
where  
\be G =\left[ (r+1) \eps_g + \gamma + \frac{ {\nu^2}  L \gamma}{(\nu - 1 )rM} \right] \frac{ \nu^2  \gamma}{(\nu - 1)rM}, \label{eq:Gdef} \ee
and $\gamma$ is defined in \cref{eq:gammap}--\cref{eq:deldef} with $\mu={\eps_g}/2 $.
\end{theorem}
\begin{proof}
{The proof is based on the observation that, when the iterates leave $C_1$, if the trust region is large enough, then by \Cref{noisyfuncred} the noisy objective function starts decreasing immediately (Case 1); otherwise the smallness of the trust region limits the increase in the objective function before the trust region becomes large enough to ensure descent (Case 2).}
We now state this precisely.

Suppose that the $K^{th}$ step is an exiting step, i.e., $x_{K}\in C_1$ and $x_{K+1}\notin C_1$. We let $\hat K > K+1$ be the index of the first iterate that returns to $C_1$. Such a $\hat K$ exists due to \Cref{thm:GCT1}. We will prove that all iterates $x_k$ with $k \in \{K+1,\ldots ,\hat K-1\}$ are contained in $C_2$.

Since $x_k\notin C_1$ for $k \in \{K+1,\ldots ,\hat K-1\}$, we have by \cref{eq:C1} that
\be \|g_k\| > \left(r+1\right) \eps_g + \frac\beta2 \label{eq:palm2},\ee
and we have seen in \cref{eq:palm}-\cref{eq:tr-g1} that this implies that  
\be
 \|\tilde g_k\| > r\eps_g+\ve, \qquad \ k \in \{K+1, \ldots ,\hat K-1\}.\label{eq:tr-g12}
 \ee 
Also, we know that a step was taken at iterate $K$ since $x_{K}\in C_1$ and $x_{K+1}\notin C_1$, and thus applying \Cref{acceptedstep} yields
\be \tilde f_{K+1} - \tilde f_K < r(1-c_0)\eps_f \label{eq:inired} .\ee

We divide the rest of the proof according to the size of $\Delta_{K+1}$ relative to $\bar \Delta$, which is defined in \cref{eq:dhat}, i.e., 
\be \bar \Delta = \frac{\gamma}{rM}. \label{eq:dhat+}\ee

\medskip\noindent
\textbf{Case 1:}
\textbf{Suppose} $\Delta_{K+1} \geq \bar \Delta$. By \cref{eq:tr-g12} and the fact that $\nu > 1$, the conditions of \Cref{tr-lbcor2} are satisfied and thus $\Delta_k > \frac{\gamma}{\nu rM}$, for $k= K+1, \ldots , \hat K-1$. 
We can therefore apply \Cref{noisyfuncred}, with $\mu = \eps_g/2>0$, for each iterate $k = K+1,\dots,\hat K-1$ to yield{\be \tilde f(x_{K+1}) \geq \tilde f(x_{K+2}) \geq \cdots \geq \tilde f(x_{\hat K}).\ee}
Combining this result with \cref{eq:inired} we obtain 
\be \tilde f_k \leq \tilde f_{K+1} < \tilde f_K + r(1-c_0)\eps_f ,
\quad k = K+1,..,\hat K.\ee
Since $ x_K\in C_1 $ and by \cref{eq:tr-asp1}, we conclude that for $k = K+1,\ldots,\hat K$,
\be f_k < f_K +[2+ r(1-c_0)]\eps_f \leq \sup_{y\in C_1}  f(y) +[2+ r(1-c_0)]\eps_f .\label{eq:P2case1}\ee
{Therefore, the inequality in \cref{eq:c2def} is satisfied in this case.}

\medskip\noindent
\textbf{Case 2:} \textbf{Suppose} $ \Delta_{K+1} < \bar \Delta$. We begin by considering the increase in the function value while the trust region remains less than $\bar \Delta.$ To this end, we define
\be l = \left\lceil{\log_\nu{\frac{\bar\Delta}{\Delta_{K+1}}} }\right\rceil,\quad 
\label{eq:defl}\ee
where $\lceil{\cdot}\rceil$ denotes the ceiling operation. Since the trust region radius is increased by a factor of at most $\nu$, we have that $l$ is the minimum number of steps required for the trust region radius to increase from $\Delta_{K+1}$ to (at least) $\bar\Delta$. Now, if $ K+ l > \hat K$, then the iterates return to $C_1$ before the trust region becomes at least $\hat \Delta$. Therefore, the number of out-of-$C_1$ iterations taken by the algorithm while $\Delta_k < \hat \Delta$ is
\begin{equation} \label{eq:lhat}
{\hat l = \min\{l - 1, \hat K - K -1\}.}
\end{equation}
The increase in function values for iterations indexed by $k = K+1, \ldots, K+\hat l + 1$ is bounded as follows:
 \begin{align}
| f(x_{k }) - f(x_{K})| &\leq \sum_{i = 0}^{k-K-1} | f(x_{K+1+i}) - f(x_{K+i})| \nonumber \\
							 &\leq  \sum_{i = 0}^{\lh} | f(x_{K+1+i}) - f(x_{K+i})| \nonumber\\
							 &\leq  \sum_{i = 0}^{\lh} \Delta_{K+i} \max_{x\in[x_{K+i},x_{ K+1+i}]}							  \| g(x)\| \nonumber \\
							 &=  \sum_{i = 0}^{\lh} \Delta_{K+i} \max_{x\in[x_{K+i},x_{ K+1+i}]}							  \| g(x) - g(x_{K+i}) + g(x_{K+i})\| \nonumber\\
							 &\leq  \sum_{i = 0}^{\lh} \Delta_{K+i} \big[ \| g(x_{K+i})\| + L\Delta_{K+i} \big] \qquad\mbox{(by \cref{eq:lip})}. \label{eq:funcstep1}
\end{align}
To estimate the right hand side, we need to bound the total displacement made by the algorithm during those iterations.
It follows from \cref{eq:defl} that
\be
\bar\Delta/\nu\leq \nu^{l-1} \Delta_{K+1}< \bar\Delta\leq\nu^l\Delta_{K+1}, \label{eq:tr@l}\ee
and thus for $i = 0, ..., {\hat l}$, 
\be \Delta_{K+1+i} \leq \nu^i\Delta_{K+1}\leq \nu^{\hat l}\Delta_{K+1} \leq \nu^{l-1}\Delta_{K+1}< \bar\Delta .\label{eq:cof_p}\ee
 By \cref{eq:tr-g12}, \cref{eq:cof_p}, we can apply \Cref{tr-lb} to each iterate $i = 0, ..., \hat l$, and obtain
\be \Delta_{i+1} = \nu \Delta_i. \ee
Thus for $i = 0, ..., \hat l$,
\be\Delta_{K+1+i} = \nu^i\Delta_{K+1}\leq \nu^{\hat l}\Delta_{K+1} \leq \nu^{l-1}\Delta_{K+1}< \bar\Delta . \label{eq:cof}\ee
Summing from $i = 0$ to $\lh$, we have
\be\label{eq:newdisplacement}
\sum_{i = 0}^{\lh} \Delta_{K+1+i} 
= \sum_{i = 0}^{\lh} \nu^{i}\Delta_{K+1}
< \frac{\bar\Delta}{\nu^{\hat l}} \sum_{i = 0}^{\lh} \nu^{i}
  = \frac{\bar\Delta}{\nu^{\lh}} \frac{ \nu^{\lh+1}-1}{\nu-1} 
  < \frac{\bar\Delta}{\nu^{\lh}} \frac{ \nu^{\lh+1}}{\nu-1}
  = \frac{ \nu}{\nu-1}\bar\Delta .
\ee
By assumption, $\Delta_{K+1}< \bar\Delta$, which implies $\Delta_K<\nu\bar\Delta$; adding this to \cref{eq:newdisplacement} we obtain
\be\sum_{i = 0}^{\lh + 1} \Delta_{K+i} < \frac{ \nu^2}{\nu-1}\bar\Delta. \label{eq:sumdel} \ee 
Therefore, for $i=0, \ldots, \lh$,
 \begin{align}		
\| g(x_{K+i})\| +L\Delta_{K+i} &=  \| g(x_{K}) + \sum_{j = 0 }^ {i-1}\left[ g(x_{K+j+1}) - g(x_{K+j})\right] \| +L\Delta_{K+i} \nonumber \\
		   &\leq  \left\| g(x_{K}) \right\| + \sum_{j = 0 }^ {i-1}\left\| g(x_{K+j+1}) - g(x_{K+j})\right\|+L\Delta_{K+i} \nonumber\\
		   &\leq  \left\| g(x_{K}) \right\| + \left(\sum_{j = 0 }^ {i-1}L\Delta_{K+j}\right)+L\Delta_{K+i} \nonumber \\
		   &<  \left\| g(x_{K}) \right\| + L\sum_{j = 0 }^ {\hat l + 1 }\Delta_{K+j} \quad (\mbox{since} \ i< \lh +1) \nonumber\\
		   &<  \left\| g(x_{K}) \right\| + \frac{\nu^2}{\nu - 1} L \bar \Delta \qquad (\mbox{by \cref{eq:sumdel}}). \label{eq:gradbound1}
\end{align}
Substituting this inequality into \cref{eq:funcstep1}, we obtain for any $k = K+1,..., K+\hat l + 1$,
 \begin{align}
| f(x_{k }) - f(x_{K})|    &\leq \sum_{i = 0}^{\lh} \Delta_{K+i} \left[\left\| g(x_{K}) \right\| + \frac{\nu^2}{\nu - 1}  L \bar \Delta\right] \nonumber \\
						 & <  \left[ \| g(x_{K})\| + \frac{\nu^2}{\nu - 1}  L \bar\Delta \right] \frac{\nu^2}{\nu - 1} \bar\Delta \nonumber \\
					 	 & \leq  \left[ (r+1) \eps_g + \gamma + \frac{\nu^2}{\nu - 1}  L \bar\Delta \right] \frac{\nu^2}{\nu - 1}  \bar\Delta \qquad(\mbox{since $x_K \in C_1$}) \nonumber \\
					     & = \left[ (r+1) \eps_g + \gamma + \frac{ {\nu^2}  L \gamma}{(\nu - 1 )rM} \right] \frac{ \nu^2  \gamma}{(\nu - 1)rM} \qquad\mbox{(by \cref{eq:dhat+})} \nonumber \\
						 & =  G. \label{eq:func_bd}
\end{align}
Therefore, {for $k=K+1, \ldots, K+ \hat l + 1$,}
\be f(x_{k}) < f(x_{K}) + G\leq \sup_{y\in C_1} f(y) + G . \label{eq:func_bd_final}\ee
We now consider two possibilities.

\medskip\noindent
\textbf{Case 2a): Suppose} $ K+1+l > \hat K$. Then, $\hat K - K -1 \leq l-1$ and by \cref{eq:lhat} we have that $ \hat l =\hat K - K -1 $. Condition \cref{eq:func_bd_final}, thus reads 
\be f(x_{k}) < f(x_{K}) + G\leq \sup_{y\in C_1} f(y) + G, \qquad
k=K+1, \ldots,\hat K  , \label{eq:before2}\ee
and thus the inequality in \cref{eq:c2def} is satisfied for $k=K+1, \ldots,\hat K -1 $.

\medskip\noindent
\textbf{Case 2b): suppose} $K+1+l \leq \hat K$. Then, by \cref{eq:defl} we have that $ \hat l = l - 1$, and \cref{eq:func_bd_final} reads
\be f(x_{k}) < f(x_{K}) + G\leq \sup_{y\in C_1} f(y) + G \qquad k=K+1, \ldots, K + l . \label{eq:before1}\ee
Let us now consider the iterates following $K+l$ that are outside $C_1$, i.e., those indexed by $ k = K+l+1, ..., \hat K -1$. Letting $i = \lh = l-1$ in \cref{eq:cof} and recalling the first inequality in \cref{eq:tr@l},
\be \Delta_{K+l} = \nu^{l-1} \Delta_{K+1} \geq 
\frac{\bar\Delta}{\nu}.\label{eq:tr@lineq}\ee 
We can therefore apply \Cref{tr-lbcor2} to iterates indexed by $ k = K+l+1, ..., \hat K -1$ and deduce that 
\[ \Delta_k \geq \frac{\bar\Delta}{\nu}, \quad k= K+l+1, ..., \hat K -1.
\] 
This fact, together with \cref{eq:tr-g12}, allow us to invoke \Cref{noisyfuncred}, for $ k = K+l, ..., \hat K -1$, to yield
\be \tilde f(x_{K+l}) \geq \tilde f(x_{K+1+l}) \geq \tilde f(x_{K+2+l}) \geq ... \geq\tilde f(x_{\hat K}).\label{eq:red_relation}\ee
Recalling \cref{eq:before1} with $k = K+l$ and using \cref{eq:tr-asp1} we obtain
\be \tilde f(x_{K+l}) < \sup_{y\in C_1} f(y) +G + \eps_f.\label{eq:nfuncl}\ee
This condition together with \cref{eq:red_relation} yields
\be f(x_{k})\leq \tilde f(x_{K+l}) + \eps_f< \sup_{y\in C_1} f(y) +G + 2\eps_f \qquad k = K+l, ..., \hat K-1 \label{eq:after1} .\ee
Combining this bound with \cref{eq:before1} we conclude 
\be f(x_{k})< \sup_{y\in C_1} f(y) +G + 2\eps_f, \qquad k = K+1,..., \hat K-1 ,\label{eq:P2case2a} \ee
and thus {the inequality in \cref{eq:c2def} is satisfied.}


\end{proof}

The constant $G$ defined in \cref{eq:Gdef} is proportional to $\eps_g^2,\eps_g\sqrt{\eps_f}, \eps_f$. Since that $G$ characterizes the function value bounds, the dependence on $\eps_f$ is expected; the dependence on $\eps_g$ and $\eps_g\sqrt{\eps_f}$ arises from the combined effect of the trust region radius and gradient norm.



\section{Numerical Experiments}
\label{sec:numerical}
To illustrate the performance of the proposed Algorithm~\ref{algorithm}, we coded it in {\sc matlab}  and applied it to a small selection of unconstrained optimization problems. We injected uniformly distributed noise in the evaluations of the function and gradient.  
Specifically, we let (c.f. \cref{eq:tr-noise}) 
\begin{equation} \label{eq:noiseu} 
\delta_f = X_f\in \mathbb{R}, \quad X_f \sim U(-\eps_f,\eps_f),\quad\mbox{and}\quad
 \delta_g = X_g\in \mathbb{R}^{n}, \quad X_g\sim \mathbb{B}_n(0,\eps_g), 
\end{equation} 
where $U(-a,a)$ denotes the uniform distribution from $-a$ to $a$, and $\mathbb{B}_n(0,a)$ denotes the $n$ dimensional ball centered at $0$ with radius $a$. By generating noise in this way we satisfy \Cref{assum2}. 

We set the parameters in Algorithm~\ref{algorithm} as follows: $c_0 = 0.1, c_1 = 1/4, c_2 = 1/2$ and $\nu = 2$. The solution of the trust region subproblem (Step~3 of Algorithm~\ref{algorithm}) was computed using the standard Newton-CG method described e.g. in \cite{mybook}, with termination accuracy $10^{-8}$. In order to better illustrate the performance of the algorithm in the presence of noise, we did not include a stop test and simply ran it for 200 iterations, which was sufficient to observe its asymptotic behavior.

\subsection{Failure of the Classical Trust Region Algorithm}
\begin{figure}[htp!]
\begin{center}
\includegraphics[width=.9\textwidth]{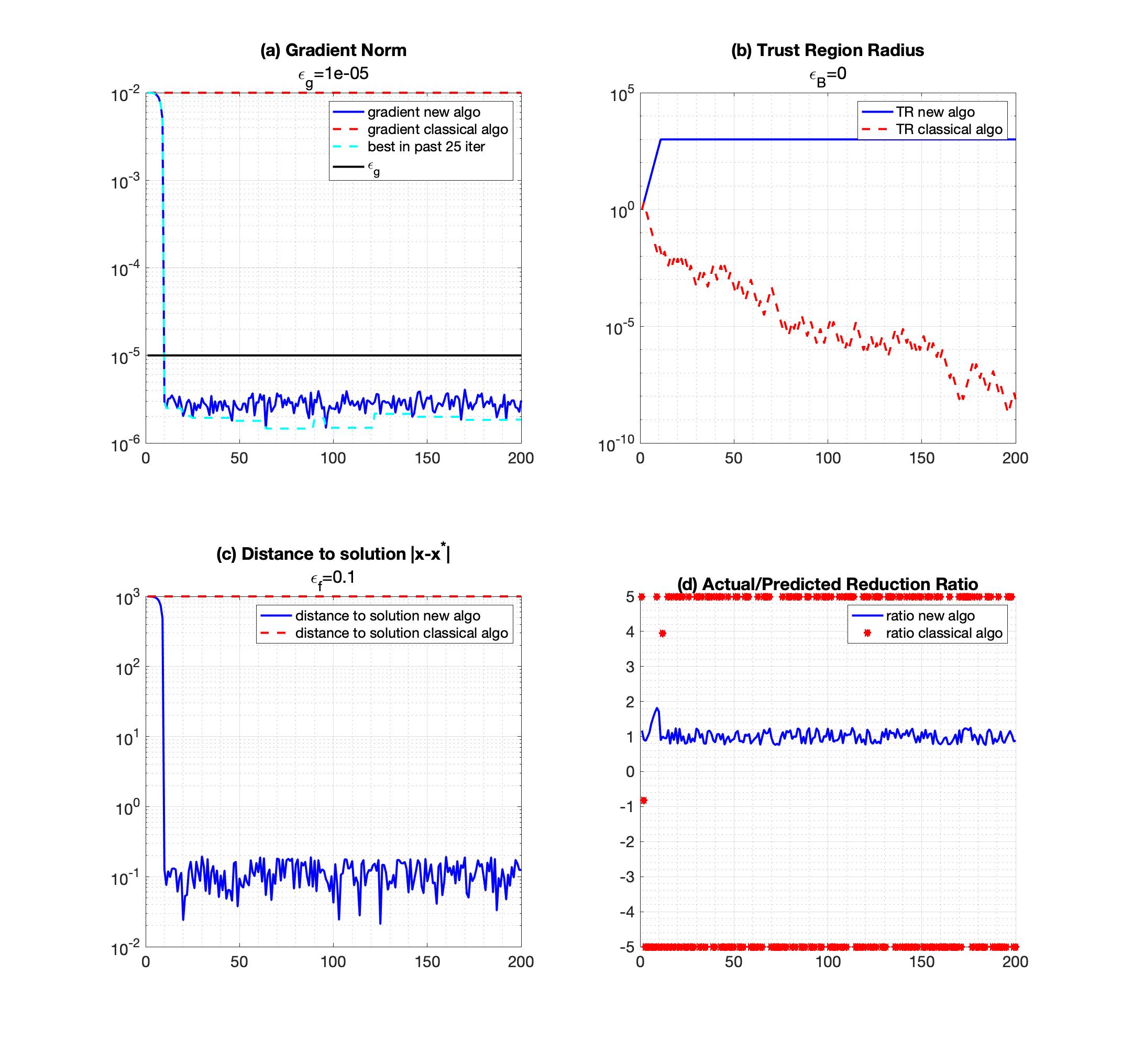}
\caption{New and classical trust region algorithms applied to a simple quadratic problem.}
\label{fig:exp_quad}
\end{center}
\end{figure}

We present two examples showing failure of the classical trust region algorithm, in contrast with Algorithm~\ref{algorithm}.
First, we consider the simple quadratic function 
\be f = x^T D x, \ee
where $x\in \cR^8$ and $D$ is the diagonal matrix
\be D = {\rm diag} (1e-5, 1e-4.75, 1e-4.5,....,1e-3.25). \ee
The condition number of $D$ is roughly $56$. We set $\eps_f = 10^{-1}$ and $\eps_g = 10^{-5}$ in \cref{eq:noiseu}. The Hessian of the quadratic model \cref{eq:tr-nmodel} was defined as $B_k = \nabla^2 f(x_k)$; i.e., we did not inject noise in this experiment. We started both algorithms from $x_0 = (1000, 0,0,....,0)$, with an initial trust region radius $\Delta_0=1$. 
The results are displayed Figure~\ref{fig:exp_quad}. 

The four panels in Figure~\ref{fig:exp_quad} compare the performance of the classical algorithm (red dashed line) and Algorithm~\ref{algorithm} (blue solid line). The horizontal axis in each panel records the iteration number. In the upper left panel (a) we report the norm of the (noiseless) gradient $\| \nabla f(x_k)\|$, along with the injected noise level $\eps_g$ (solid black line); the light blue dashed line plots the lowest value generated by Algorithm~\ref{algorithm} in the past 25 iterations. In the upper right panel (b) we report the trust region radius; in the lower-left panel (c) the distance to solution; and in the lower right panel (d), the computed actual-to-predicted reduction ratio $\rho_k$; for graphical clarity, ratios greater than $5$ or less than $-5$ were plotted as $+/- 5$ in panel (d).

We observe that the classical algorithm exhibits large oscillations in $\rho_k$, which causes the trust region radius to shrink so much that significant progress cannot be made. In contrast, $\rho_k$ is controlled well in Algorithm~\ref{algorithm}. In this test, initial the trust region radius $\Delta_0$ is not small.

\begin{figure}[h!]
\begin{center}
\includegraphics[width=.9\textwidth]{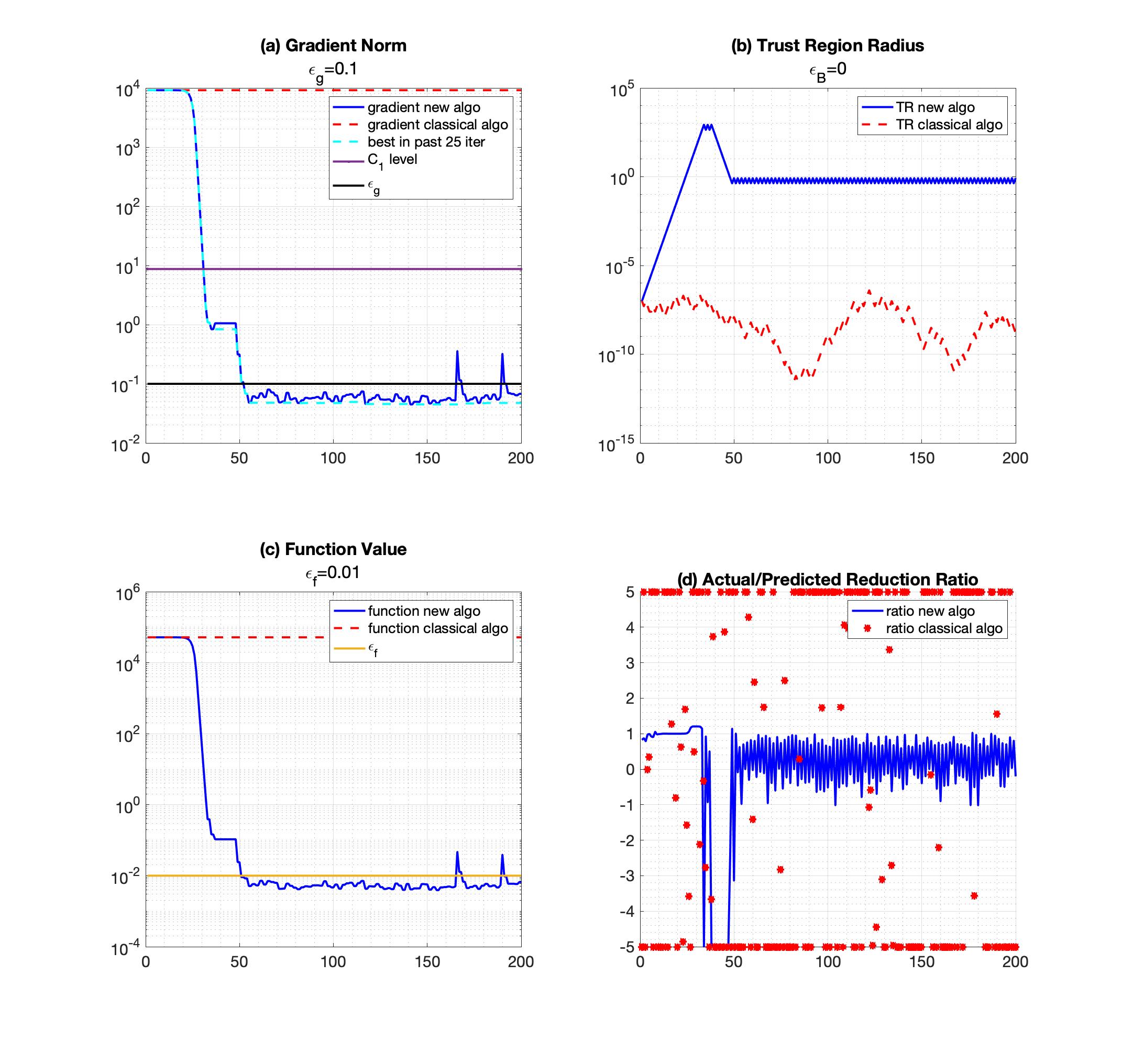}
\caption{New and classical trust region algorithms initialized with small trust region radius.}
\label{fig:exp_sanity_check}
\end{center}
\end{figure}

In the next experiment, we illustrate the damaging effect that a very small $\Delta_0$ can have on the classical algorithm, but not on the proposed algorithm. We applied the two algorithms to the following tri-diagonal function
\be
f(x)= \frac{1}{2}\left(x^{(1)}-1\right)^{2}+\frac{1}{2} \sum_{i=1}^{N-1}\left(x^{(i)}-2 x^{(i+1)}\right)^{4}, \quad N=200.\label{eq:prob}
\ee
The results are reported in Figure~\ref{fig:exp_sanity_check}. In the upper left panel, we additionally plot in purple the size of the critical region $C_1$, i.e. the value of the right-hand side in \cref{eq:C1}. (The latter requires knowledge of the constant $M$, which we approximate by the norm of the Hessian at the solution.) This panel shows that the theoretical prediction given in \Cref{thm:GCT1} is pessimistic when compared to the final achieved accuracy in the gradient, as is to be expected of convergence results that assume that the largest possible error occurs at every iteration.
The upper right hand panel illustrates that Algorithm~\ref{algorithm} is able to quickly increase the trust region radius an allow progress, unlike the classical algorithm.

\subsection{General Performance of the Proposed Algorithm}
We also tested the two algorithms on a subset of problems from \cite{schittkowski1987more}; the results are presented in the {supplementary material}. As a representative of these runs, we report the results for the tri-diagonal objective function \cref{eq:prob}.
This time, the Hessian $ B_k$ of the quadratic model \cref{eq:tr-nmodel} is obtained by injecting noise in the true Hessian matrix. We define
\be \label{eq:noiseb}  B_k = \nabla^2 f(x_k) + \delta_B,\ee
\be \delta_B = \frac{A^T \Lambda A}{\|A\|^2}, \quad A_{ij}\sim U(0,1), \quad (\Lambda)_{ii} \sim U(-\eps_B,\eps_B),\ee
{where $\Lambda$ is a diagonal matrix.}
Thus, the matrices $B_k$ are symmetric but not necessarily positive definite.
We employed larger noise levels than in the previous experiments: $\eps_f = 10$, $\eps_g = 100$, and $\eps_B = 1000$. This simulates the situation that may occur when employing finite difference approximations, where the error increases with the order of differentiation.
Both algorithms were initialized from the same starting point $x_0$, which was generated such that each entry in $x_0$ is sampled uniformly from $-50$ to $50$. To ensure a fair comparison, at each iterate we inject exactly the same noise into both algorithms.  

\begin{figure}[h!]  
\begin{center}
\includegraphics[width=.9\textwidth]{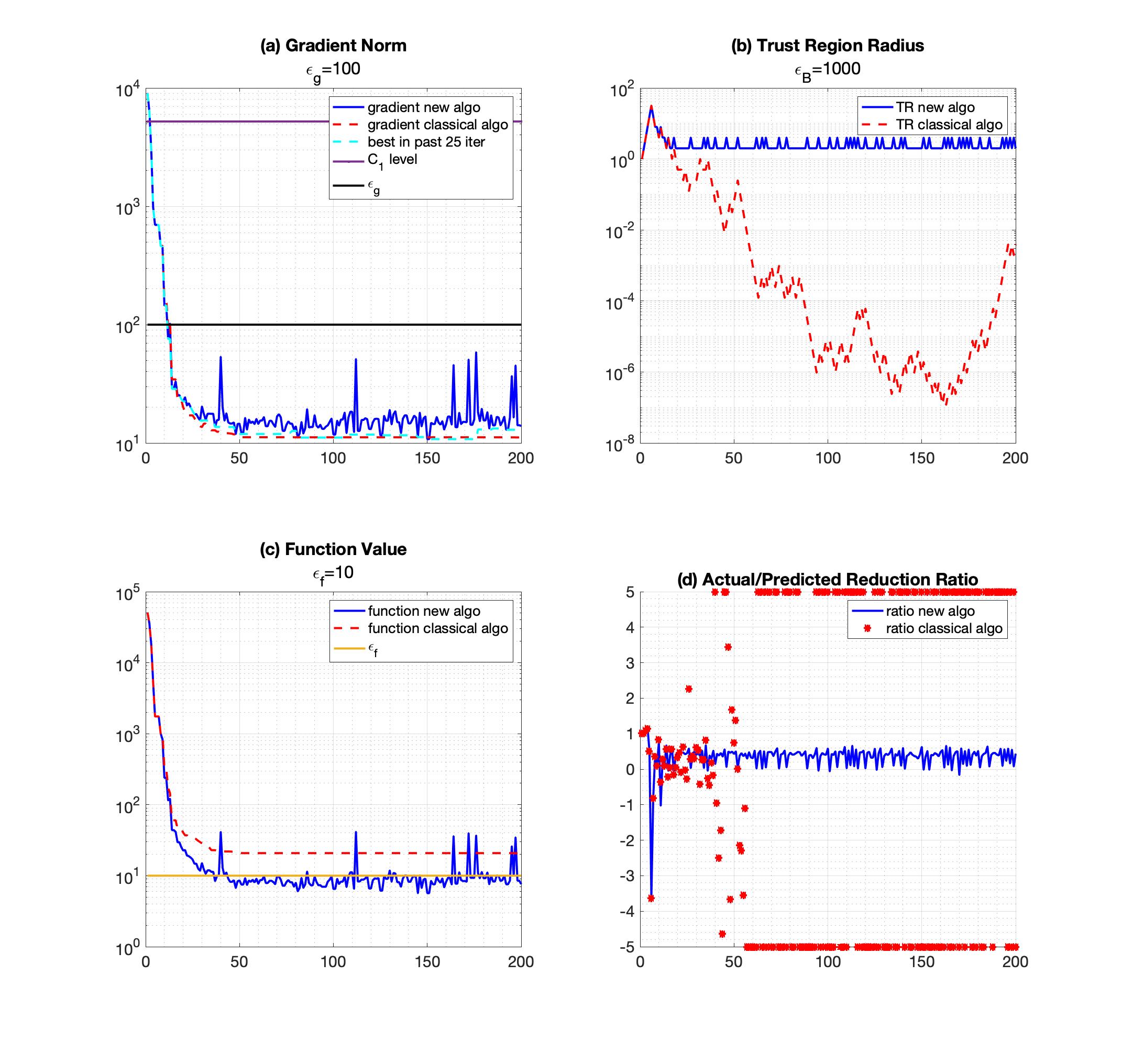}
\caption{Comparison of the new and classical trust region algorithms when solving problem \cref{eq:prob} with uniform noise given by \cref{eq:noiseu} \cref{eq:noiseb}. }
\label{fig:exp_uniform}
\end{center}
\end{figure}

We report the results in Figure \ref{fig:exp_uniform}, which displays the same information as in Figure~\ref{fig:exp_sanity_check}. We observe that both algorithms perform similarly before entering the noisy regime. Algorithm~\ref{algorithm} exhibits larger oscillations in the gradient norm due to the larger trust region radius, but achieves a lower objective function value. Whereas the large reduction in the trust region radius led to failures of the classical algorithm in the examples reported above, in many test runs such as that given in Figure~\ref{fig:exp_uniform}, it can be beneficial by producing increasingly smaller steps that yield milder oscillations in the gradient norm than Algorithm~\ref{algorithm}. We cannot, however, recommend this type of trust region reduction as a general procedure for handling noise since failures can happen unexpectedly. 

\subsection{Evaluating the Theoretical Results}\label{sec:epsf}

We have seen that the critical region $C_1$ gives a pessimistic estimate of the achievable accuracy in the gradient because the analysis assumes worst-case behavior at each iteration, rather than providing estimates in high probability. Nevertheless, \Cref{thm:GCT1} identifies the functional relationship between the achievable accuracy and the noise level: the right hand side in \cref{eq:C1} scales as a function of $\eps_g$ and $\sqrt{\eps_f}$. We performed numerical tests to measure if the accuracy achieved in practice scales in that manner.

\begin{figure}[h!]
\begin{center}
\includegraphics[width=\textwidth]{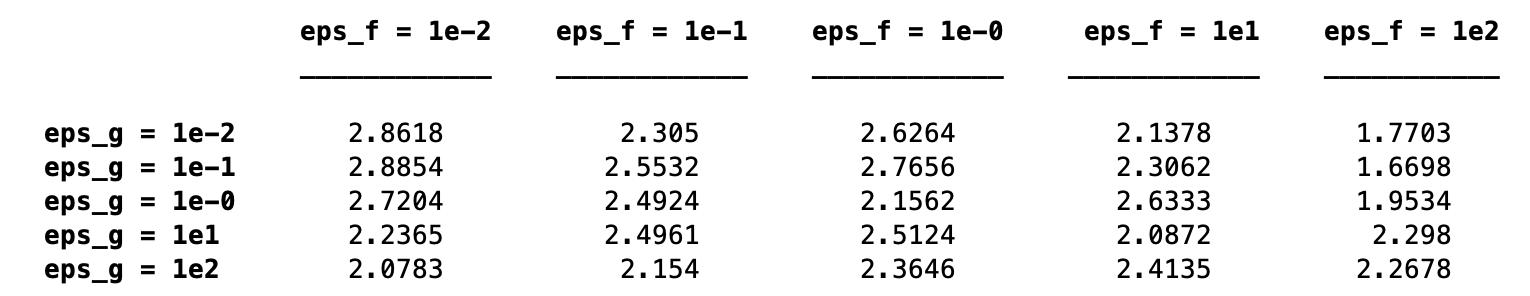}
\caption{$R(\eps_f,\eps_g)$ given in \cref{eq:rr}: Log$_{10}$ of the ratio between predicted and actual accuracy in the gradient, as a function of these noise level $\epsilon_f, \epsilon_g$. 
The small variation in these numbers suggests that \Cref{thm:GCT1} gives the correct dependence on the noise levels. }
\label{fig:table}
\end{center}
\end{figure}

We employed the tridiagonal function \cref{eq:prob}, for which we can estimate the constant $M$, as mentioned above. For given $\eps_f$ and $\eps_g$, we compute the right hand side in \cref{eq:C1}, {which we denote as $C(\eps_f,\eps_g)$, and ran Algorithm~\ref{algorithm} as in the previous test. We repeated the run 10 times using different seeds, $s=1,\ldots,10$, to generate noise. For each run, we track the smallest value} of $\| \tilde g_k\|$ during the most recent 25 iterations and record the smallest such value observed during the run, which we denote as $\|\tilde g^*({\eps_g,\eps_f,s})\|$, where $s$ denotes the seed. In \Cref{fig:table}, we report the quantity
\be R(\eps_f,\eps_g) =\log_{10} \frac{C(\eps_f,\eps_g)}{\sum_{s = 1}^{10} \|\tilde g^*({\eps_g,\eps_f,s})\|} \label{eq:rr} \ee
as we vary $\eps_f$ and $\eps_g$ from $10^{-2}$ to $10^2$. 
The fact that the ratio between the theoretical bound and the smallest gradient norm measured in practice remained roughly constant gives numerical support to the claim that the achievable gradient norm is proportional to $\eps_g$ and $\sqrt{\eps_f}$. We should note that these observations are valid only when averaging multiple runs with different seeds, as one can observe significant variations among individual runs of Algorithm~\ref{algorithm}.

\section{Final Remarks}
\label{sec:conclusions}

In this paper, we proposed a noise-tolerant trust region algorithm that avoids the pitfall of the classical algorithm, which can shrink the trust region prematurely, preventing progress toward a stationary point. Robustness is achieved by relaxing the ratio test used in the step acceptance, so as to account for errors in the function.

We showed that when the noise in the function and gradient evaluations is bounded by the constants $\epsilon_f, \epsilon_g$, an infinite subsequence of iterates satisfies
\begin{equation} \label{eq:heart}
    \| g_k \| = O({\sqrt{\epsilon_f}}, \epsilon_g). \end{equation}
When noise is not present, our results yield the limit $\{ \|g_k\|\} \rightarrow 0$ (the sets $C_1$ and $C_2$ in \Cref{thm:GCT1} and \Cref{thm:GCT2i} coincide in this case). 

The technique and analysis presented here are relevant to the case when noise can be diminished as needed, as assumed e.g. in \cite{katya2018storm,cartis2019TRsupermartingale,bollapragada2018adaptive}. Algorithm~\ref{algorithm} can be run until it ceases to make significant progress, at which point the accuracy in the function and gradient is increased (i.e., $\epsilon_f, \epsilon_g$ are reduced) and the algorithm is restarted with the new value of $\eps_f$ in \cref{eq:tr-ratio}; this process can then be repeated. This provides a disciplined approach for achieving high accuracy in the solution using a noise-tolerant trust region algorithm.

\section*{Acknowledgments}
We thank Richard Byrd for many valuable suggestions and insights, as well as Figen \"Oztoprack, Andreas W\"achter and Melody Xuan for proofreading the paper and providing valuable feedback.

\bibliographystyle{spmpsci}
\bibliography{../../References/references}

\end{document}


\title{Supplementary Material: A Trust Region Method for the Optimization of Noisy Functions\thanks{Sun was supported by NSF grant DMS-1620022. Nocedal was supported by AFOSR grant FA95502110084 and  ONR grant N00014-21-1-2675.}
}


\author{Shigeng Sun         \and
        Jorge Nocedal 
}


\institute{Shigeng Sun \at
              Department of Engineering Sciences and Applied Mathematics, Northwestern University, Evanston, IL, USA \\
              \email{shigengsun2024@u.northwestern.edu}           
           \and
           Jorge Nocedal \at
              Department of Industrial Engineering and Management Sciences, Northwestern University, Evanston, IL, USA\\
              \email{j-nocedal@northwestern.edu}
}

\date{Received: 02 January 2022 / Accepted: }

\maketitle

\section{Additional Numerical Experiments}
{We present supplementary results on the performance of Algorithm~1.}

\subsection{Tridiagonal Function with Radamacher Noise}
\begin{figure}[h!]
\begin{center}
\includegraphics[width=\textwidth]{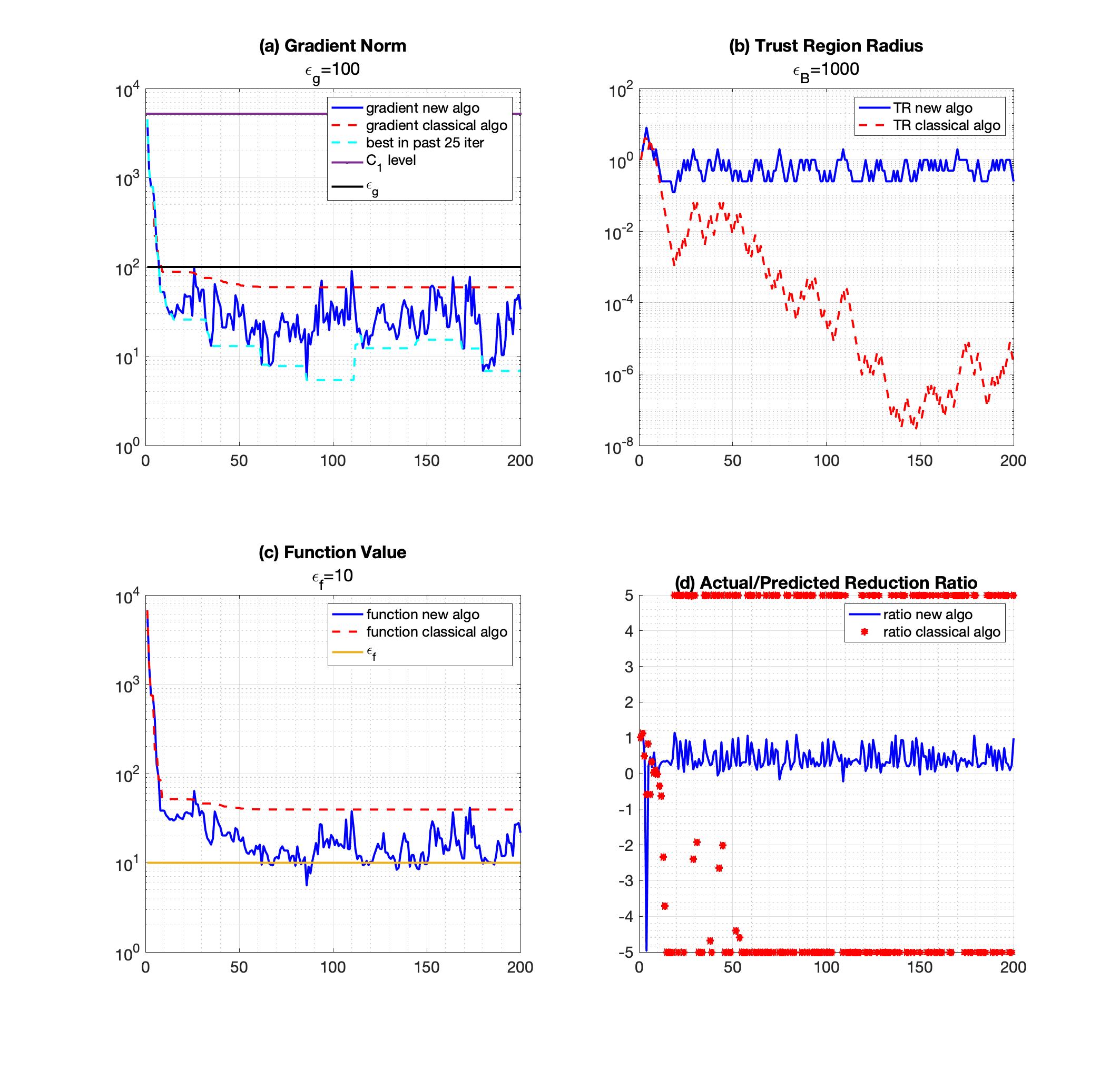}
\caption{Comparison of the proposed and classical region algorithms in the presence of Radamacher noise.}
\label{fig:exp_rad}
\end{center}
\end{figure}

In \Figref{exp_rad}, we report results of Algorithm~1 applied to the tridiagonal function described in the main paper with injected noise following the Radamacher distribution (in place of (79), (83), (84) from the main paper):
\begin{align*} 
\delta_f  =  & X_f\in\cR, \quad X_f \sim R(-\eps_f,\eps_f) \\
 \delta_g = & X_g\in\cR^{N},\quad X_g\sim \partial B_N(0,\eps_g)\\
 \delta_B = & \frac{A^T \Lambda_B A}{\|A\|^2}, \quad A_{ij}\sim  U(0,1),\quad (\Lambda_B)_{ii} \sim R(-\eps_B,\eps_B).
\end{align*}
Here $X\sim R(-a,a)$ means that the only possible values of $X$ are $\{+a, -a\}$, each with probability 1/2.

\newpage
\subsection{Tridiagonal Function with Uniform Noise}
\begin{figure}[h!]
\begin{center}
\includegraphics[width=\textwidth]{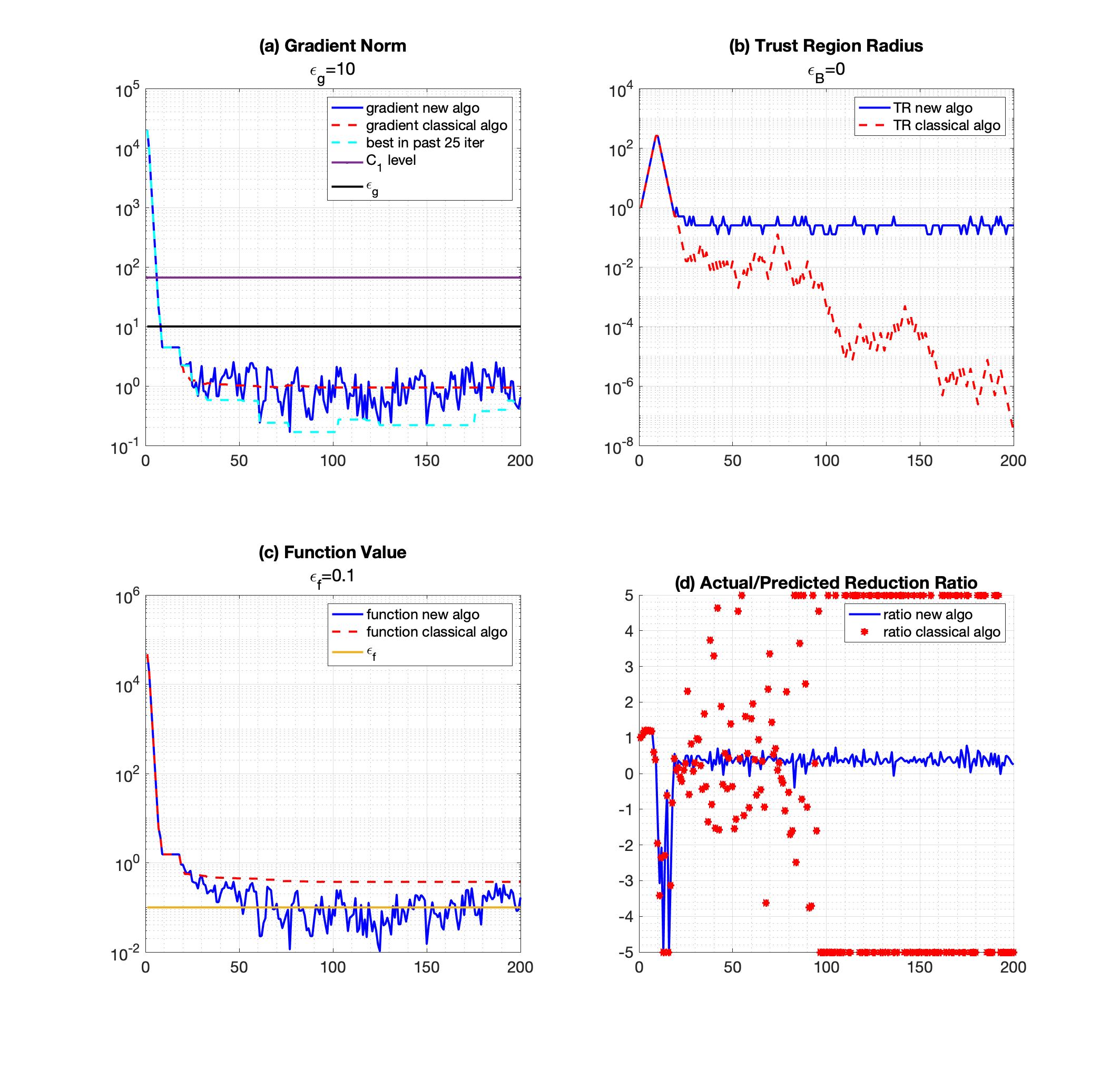}
\caption{Comparison of the proposed and classical region algorithms on the tridiagonal function with uniform noise.}
\label{fig:exp_uniform}
\end{center}
\end{figure}

\begin{figure}[h!]
\begin{center}
\includegraphics[width=\textwidth]{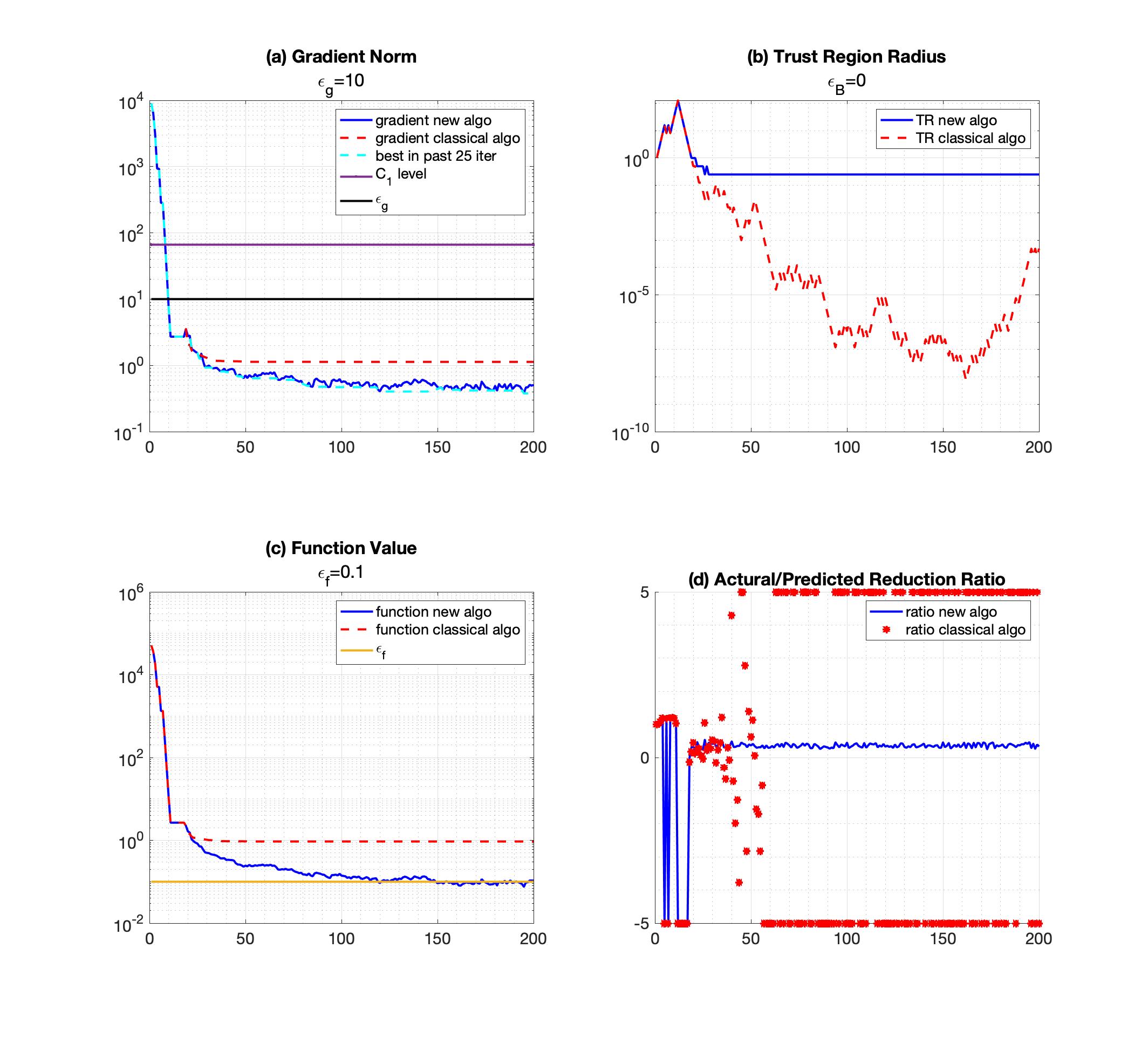}
\caption{Comparison of the proposed and classical region algorithms on the tridiagonal function with uniform noise.}
\label{fig:exp_uniform1}
\end{center}
\end{figure}

In \Figref{exp_uniform,exp_uniform1} we report some additional runs of Algorithm~1 on the tridiagonal function with different levels of uniformly distributed noise. The noise levels are given in the headers of each panel.

\newpage
\subsection{Additional Functions from Schittkowski Test Set \cite{schittkowski1987more}}
In this section, we report some additional runs on other selected problems in \cite{schittkowski1987more}. We employed the starting points given in that test set. In the following experiments, we injected  uniformly distributed noise (c.f. (79) from the main paper).

\subsubsection{Problem 271, SUR-T1-12}
\begin{figure}[h!]
\begin{center}
\includegraphics[width=\textwidth]{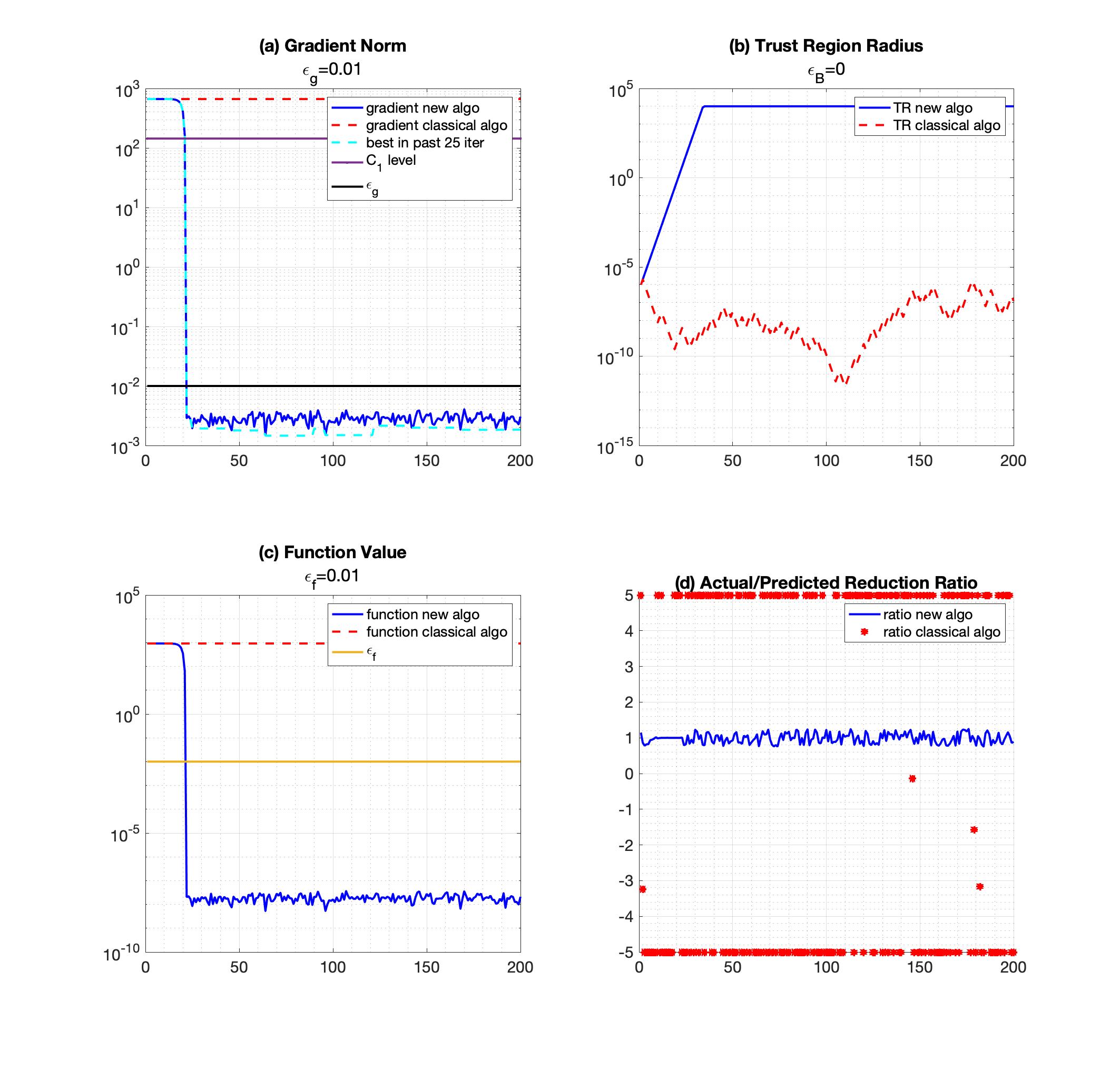}
\caption{Comparison of the proposed and classical trust region algorithms on problem 271, with a small initial trust region radius.}
\label{fig:exp_s271small}
\end{center}
\end{figure}

\begin{figure}[h!]
\begin{center}
\includegraphics[width=\textwidth]{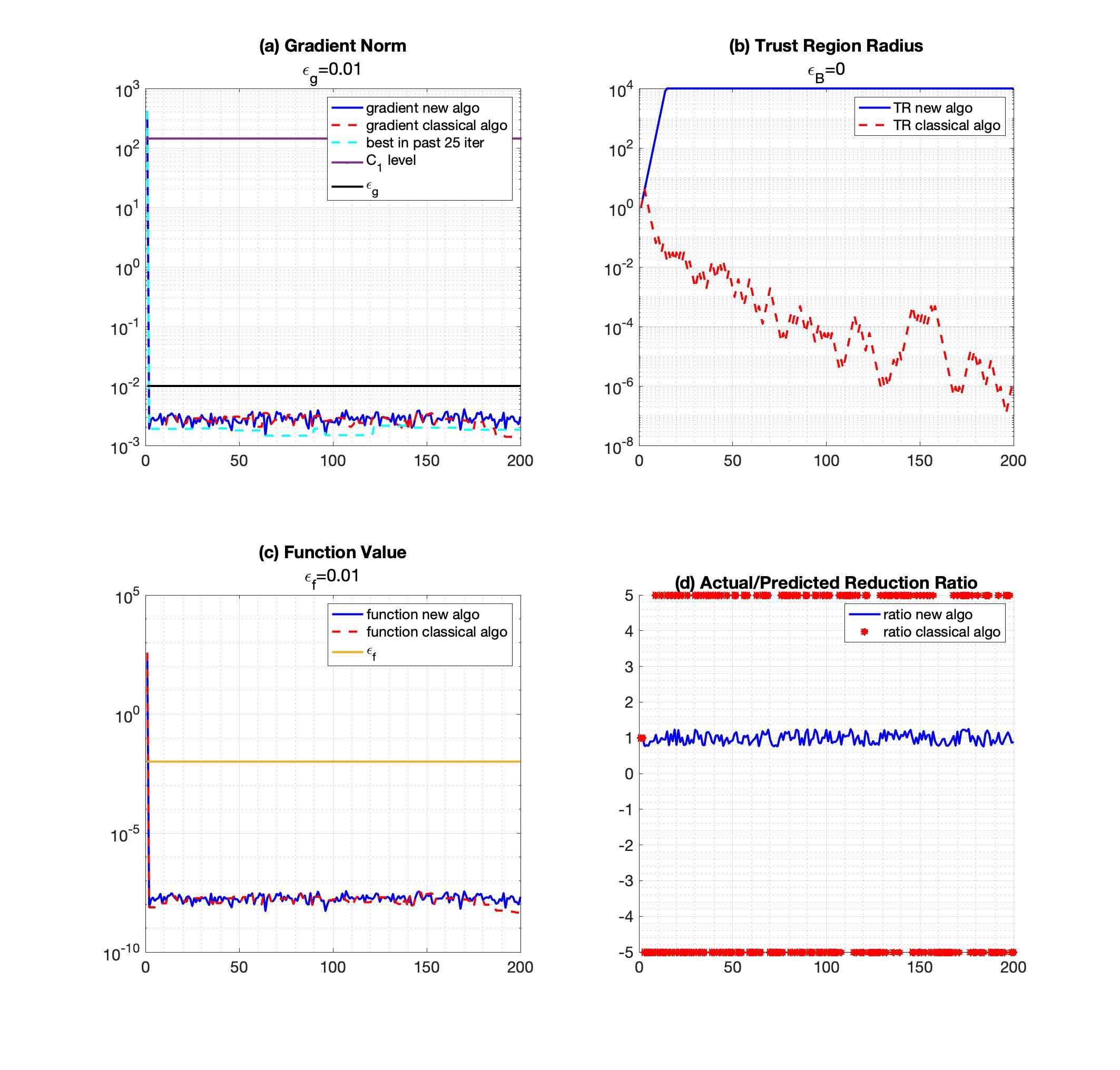}
\caption{Comparison of the proposed and classical trust region algorithms on problem 271, with a large initial trust region radius.}
\label{fig:exp_s271big}
\end{center}
\end{figure}

We started both algorithms with a small ($\Delta_0 = 1e-6$) or a large ($\Delta_0 = 1$) trust region radius, and plotted the results in \Figref{exp_s271small,exp_s271big}.

\newpage

\subsubsection{Problem 289, GUR-T1-3}
\begin{figure}[h!]
\begin{center}
\includegraphics[width=\textwidth]{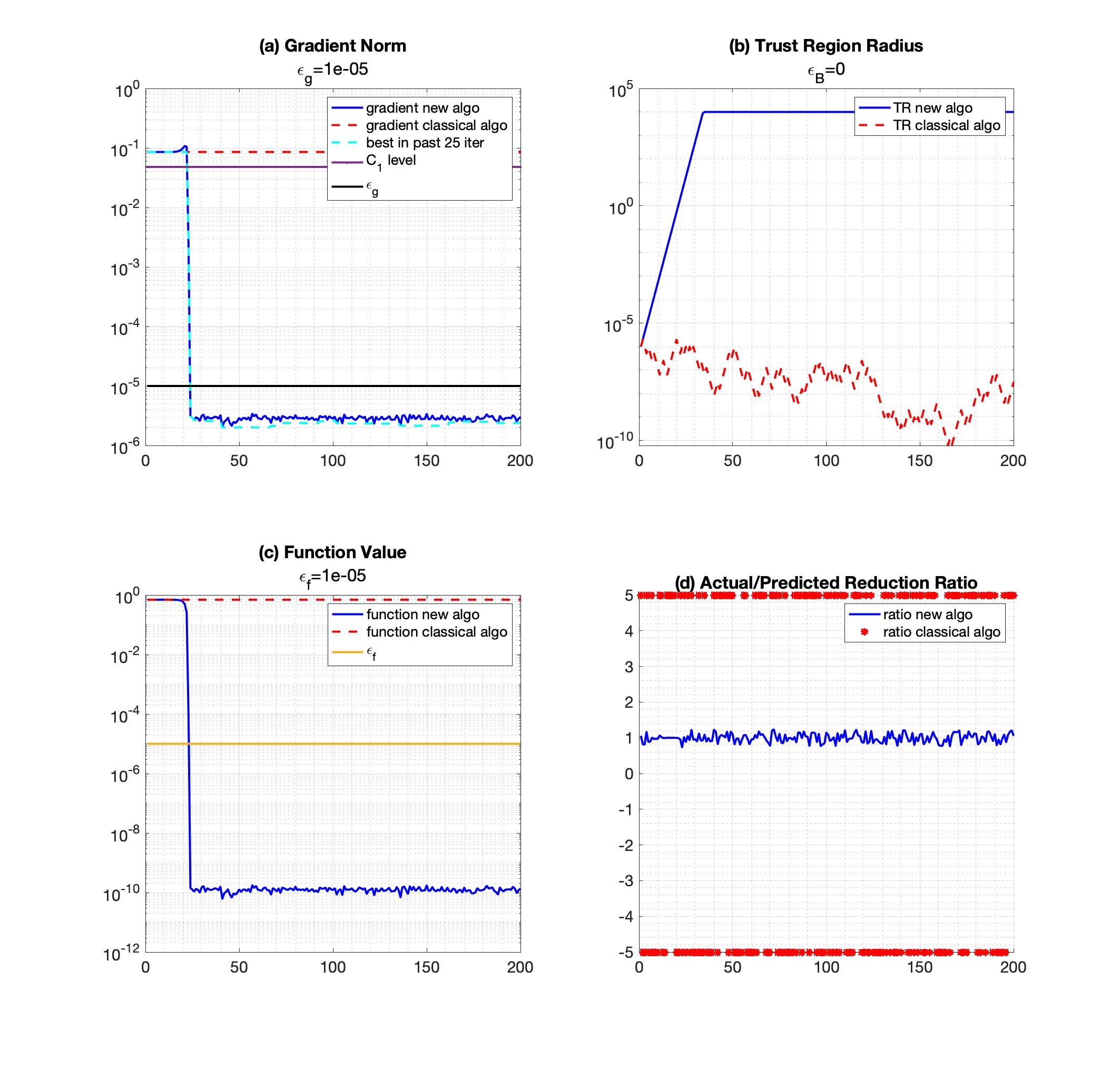}
\caption{Comparison of the proposed and classical trust region algorithms on problem 289, with small initial trust region radius.}
\label{fig:exp_s289small}
\end{center}
\end{figure}

\begin{figure}[h!]
\begin{center}
\includegraphics[width=\textwidth]{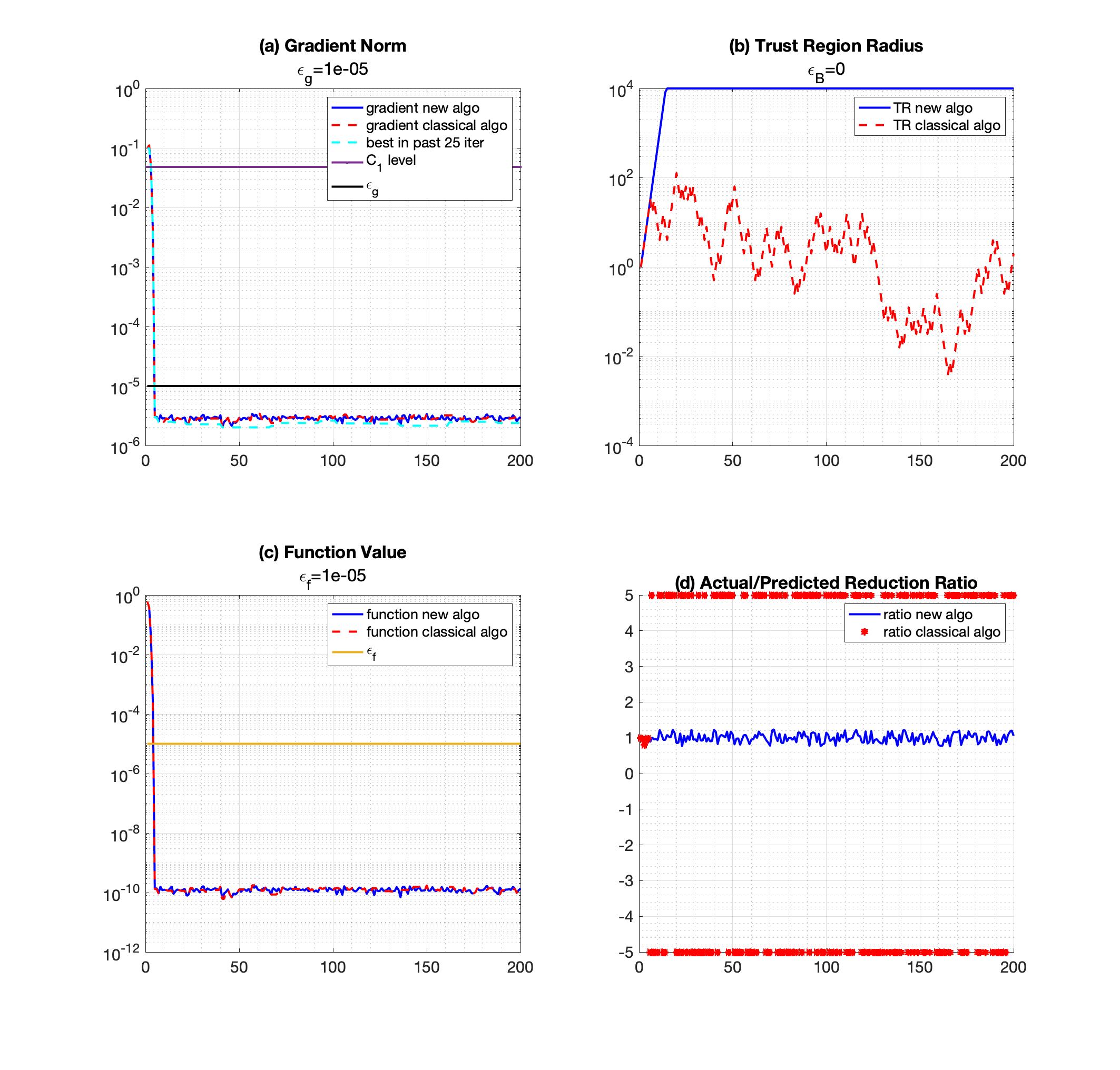}
\caption{Comparison of the proposed and classical trust region algorithms on problem 289, with large initial trust region radius.}
\label{fig:exp_s289big}
\end{center}
\end{figure}

We initiated both algorithms with small ($\Delta_0 = 1e-6$) and large ($\Delta_0 = 1$) trust region radius and plotted the results in \Figref{exp_s289small,exp_s289big}.

\newpage

\subsubsection{Problem 293, PUR-T1-18}
\begin{figure}[h!]
\begin{center}
\includegraphics[width=\textwidth]{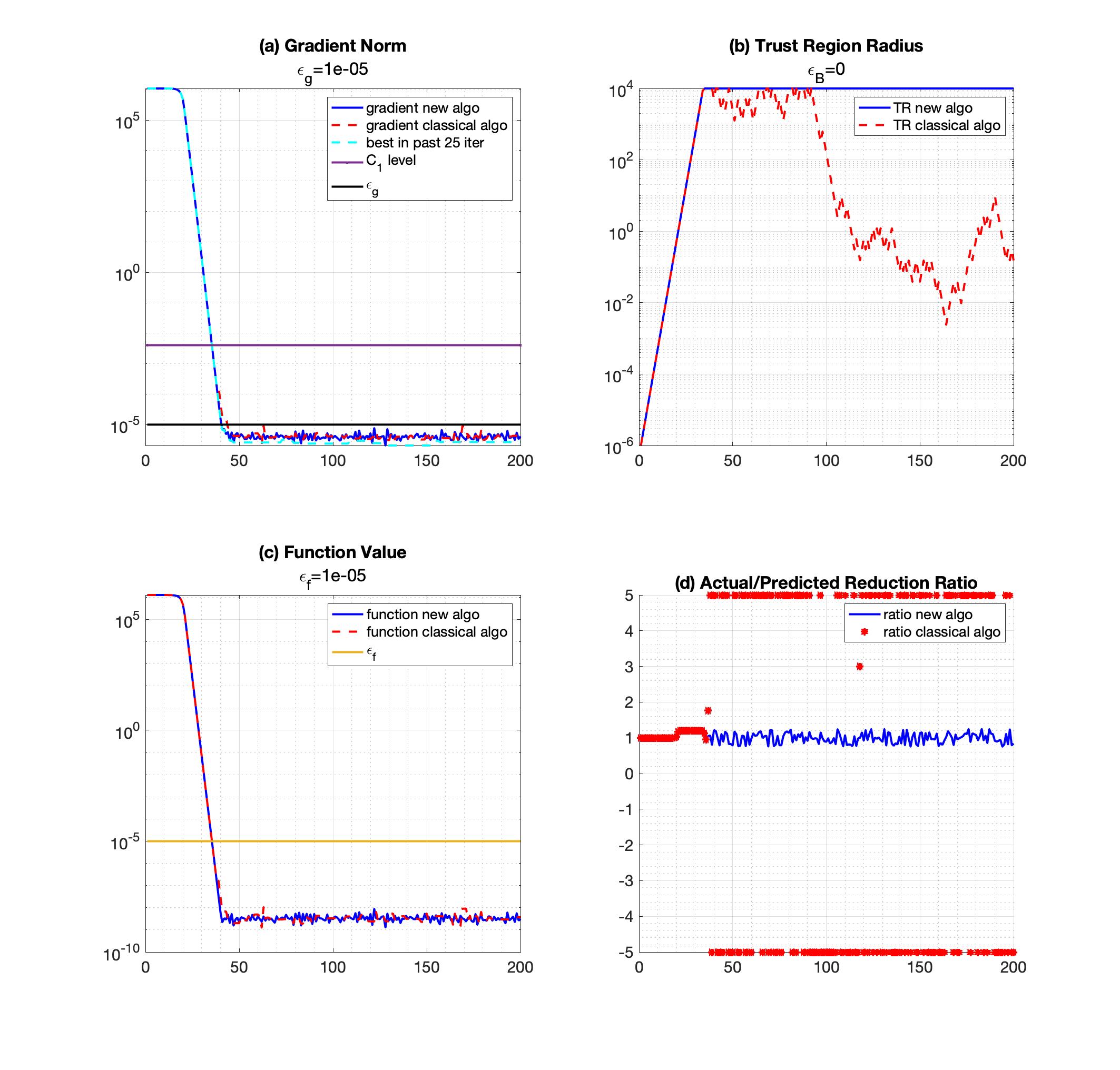}
\caption{Comparison of the proposed and classical trust region algorithms on problem 293, with small initial trust region radius.}
\label{fig:exp_s293small}
\end{center}
\end{figure}

\begin{figure}[h!]
\begin{center}
\includegraphics[width=\textwidth]{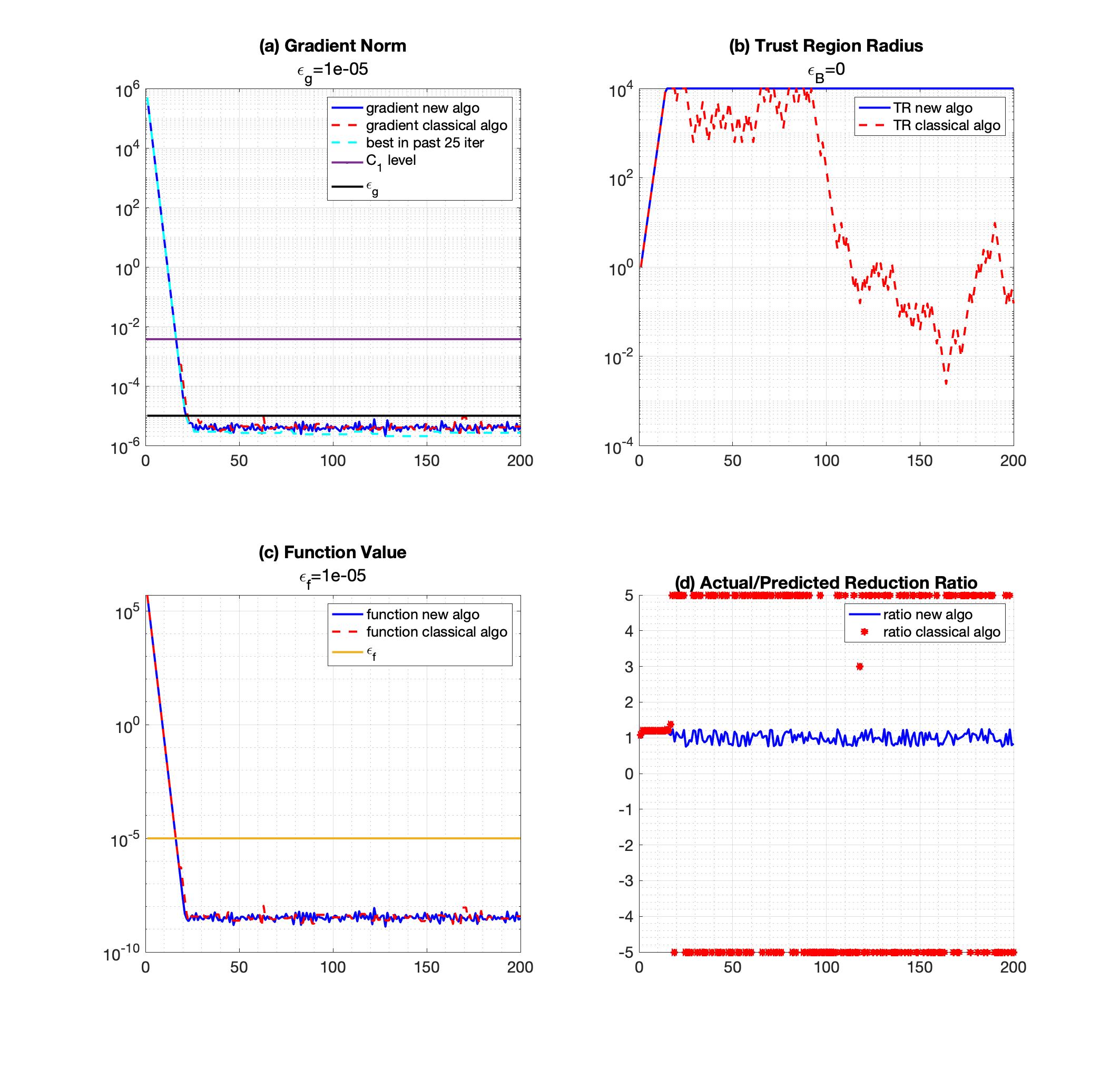}
\caption{Comparison of the proposed and classical trust region algorithms on problem 293, with large initial trust region radius.}
\label{fig:exp_s293big}
\end{center}
\end{figure}

We initiated both algorithms with small ($\Delta_0 = 1e-6$) and large ($\Delta_0 = 1$) trust region radius and plotted the results in \Figref{exp_s293small,exp_s293big}.

\newpage
\bibliographystyle{spmpsci}
\bibliography{../../References/references}